\documentclass[letterpaper, 10 pt, conference, two column]{ieeeconf}  % Comment this line out if you need a4paper

\usepackage{amsmath,amssymb,amsfonts}
\usepackage{mathrsfs}
\usepackage{setspace}
\usepackage{cite}
\usepackage{graphicx}
\usepackage{textcomp}
\usepackage{xcolor}
\usepackage{lipsum}
\usepackage{arydshln}
\usepackage{mwe}
\def\BibTeX{{\rm B\kern-.05em{\sc i\kern-.025em b}\kern-.08em
    T\kern-.1667em\lower.7ex\hbox{E}\kern-.125emX}}
    
\usepackage{mathtools}
\usepackage{cuted}
\usepackage{dsfont}
\usepackage{xspace}

\usepackage{optidef}

\usepackage{pgf,tikz}%milad
\usetikzlibrary{positioning}
\usepackage{pgfplots}
\pgfplotsset{compat = newest}
\pgfplotsset{every x tick label/.append style={font=\footnotesize, yshift=0.5ex}}
\pgfplotsset{every y tick label/.append style={font=\footnotesize, xshift=0.5ex}}
\usepackage{tikz-cd}

\newcommand{\prox}{\mathbf{prox}}

\newcommand{\sign}{\mathop{\bf sign}}

\newcommand{\ie}{{\it i.e.\;}}

\DeclarePairedDelimiter\parentheses{\lparen}{\rparen}
\newcommand{\vect}[1]{\operatorname{vec} \parentheses*{#1}}

\newcommand{\trace}[1]{\operatorname{trace} \parentheses*{#1}}

\newcommand{\ip}[1]{\left\langle #1 \right\rangle}

\newcommand{\target}{\mathcal{T}}

\tikzstyle{neuron} = [draw, fill=white, circle, 
    minimum height=3.5em, minimum width=1em, style=thick]

\newtheorem{theorem}{Theorem}
\newtheorem{proposition}{Proposition}
\newtheorem{assumption}{Assumption}
\newtheorem{corollary}{Corollary}

\newtheorem{definition}{Definition}
\newtheorem{lemma}{Lemma}

\allowdisplaybreaks

\definecolor{darkgreen}{HTML}{008000}
\definecolor{darkred}{HTML}{A52A2A}

\newcommand{\comment}[1]{}
\newcommand{\edo}{\end{document}}

\newcommand{\R}{{\mathbb R}}  %ams bold
\newcommand{\be}[1]{\begin{equation}\label{#1}}
\newcommand{\ee}{\end{equation}}
\newcommand{\halmos}{\rule{1ex}{1.4ex}}
\newcommand{\qed}{\hfill \halmos} %put \qed at right margin
\newcommand{\beq}{\begin{eqnarray}}
\newcommand{\eeq}{\end{eqnarray}}
\newcommand{\beqn}{\begin{eqnarray*}}
\newcommand{\eeqn}{\end{eqnarray*}}
\newcommand{\bi}{\begin{itemize}}
\newcommand{\ei}{\end{itemize}}
\newcommand{\ben}{\begin{enumerate}}
\newcommand{\een}{\end{enumerate}}
\newcommand{\bes}[1]{\begin{subequations}\label{#1}\begin{eqnarray}}
\newcommand{\ees}[1]{\end{eqnarray}\end{subequations}}

\DeclareFontFamily{U}{mathx}{\hyphenchar\font45}
\DeclareFontShape{U}{mathx}{m}{n}{
      <5> <6> <7> <8> <9> <10>
      <10.95> <12> <14.4> <17.28> <20.74> <24.88>
      mathx10
      }{}
\DeclareSymbolFont{mathx}{U}{mathx}{m}{n}
\DeclareFontSubstitution{U}{mathx}{m}{n}
\DeclareMathAccent{\widebar}{0}{mathx}{"73}

\newcommand{\kd}{k_d}

\newcommand{\kdst}{k_d^{*}}

\newcommand{\x}{x}
\newcommand{\SSp}{\mathcal{X}}
\newcommand{\cost}[1][\empty]{\ifx#1\empty f \else f(#1) \fi}
\newcommand{\re}[1][\empty]{\mathbb{R}^{#1}}
\newcommand{\minCost}{\underline f}
\newcommand{\CPts}{\mathcal{Z}}
\newcommand{\gradCost}[1][\empty]{\nabla \cost[#1]}
\newcommand{\sol}[1][\empty]{\ifx#1\empty \phi \else \phi(#1)\fi}
\newcommand{\xs}{\x^*}

\newcommand{\Kstbl}{\mathcal{G}}
\newcommand{\PD}[1][]{\mathbb{S}^{#1}_{++}}
\newcommand{\PSD}[1][]{\mathbb{S}^{#1}_{+}}

\newcommand{\Ksstbl}[1][]{\mathcal{G}_{#1}}
\newcommand{\SSd}{\mathcal{S}}
\newcommand{\bdKstbl}{\partial\Kstbl}
\newcommand{\e}{\epsilon}

\newcommand{\lm}{\overline\lambda}

\newcommand{\Ke}{\widetilde{K}}
\newcommand{\Abs}{\bar A^*}
\newcommand{\Ykb}[1][\empty]{\bar Y_{{\rho}#1}}
\newcommand{\I}[1][\empty]{\mathcal{I}_{\rho#1}}
\newcommand{\Pm}[1][\empty]{\mathcal{P}^{#1}}
\newcommand{\vA}[1][\empty]{\mathscr{A}^{#1}}
\newcommand{\vL}[1][\empty]{\mathscr{L}^{#1}}
\newcommand{\p}{p}
\newcommand{\q}{q}
\newcommand{\PrKe}{P_{\rho}}

\newcommand{\Pkb}[1][\empty]{\bar P_{\rho #1}}
\newcommand{\Rp}{\mathscr{R}_{\rho}}
\newcommand{\Op}{\mathscr{O}_{\rho}}

\newcommand{\Keb}{\overline K(\rho)}
\newcommand{\Bb}{\bar B}
\newcommand{\Kt}{\widetilde{K}}
\newcommand{\Kr}{\widetilde{K}(\rho)}
\newcommand{\rt}{\tilde\rho}

\newcommand{\Ksat}{\mathcal{K}_\mathrm{SAT}}

\newcommand{\gPLI}{\textit{g}P\L I\xspace}
\newcommand{\sgPLI}{\textit{sg}P\L I\xspace}
\newcommand{\lPLI}{\mbox{$\ell$P\L I}\xspace}

\IEEEoverridecommandlockouts                              % This command is only needed if 
\title{\LARGE \bf
Remarks on the Polyak-\L{}ojasiewicz inequality and the convergence of gradient systems
}

\author{Arthur C. B. de Oliveira$^{1}$, Leilei Cui$^{2}$, and Eduardo D. Sontag$^{1,3}$%
\thanks{The work of EDS and ACO was partially supported by grants AFOSR FA9550-22-1-0316 and ONR N00014-21-1-2431}%
\thanks{$^{1}$Department of Electrical and Computer Engineering, Northeastern University, USA {\tt\small a.castello@northeastern.edu}}%
\thanks{$^{2}$Massachusetts Institute of Technology, Massachusetts, USA {\tt\small llcui@mit.edu}}%
\thanks{$^{3}$Department of BioEngineering, Northeastern University, USA {\tt\small e.sontag@northeastern.edu}}
}

\begin{document}

\maketitle
\thispagestyle{empty}
\pagestyle{empty}

%%%%%%%%%%%%%%%%%%%%%%%%%%%%%%%%%%%%%%%%%%%%%%%%%%%%%%%%%%%%%%%%%%%%%%%%%%%%%%%%
\begin{abstract}

    This work explores generalizations of the Polyak-\L ojasiewicz inequality (P\L I) and their implications for the convergence behavior of gradient flows in optimization problems. Motivated by the continuous-time linear quadratic regulator (CT-LQR) policy optimization problem -- where only a weaker version of the P\L I is characterized in the literature -- this work shows that while weaker conditions are sufficient for global convergence to, and optimality of the set of critical points of the cost function, the ``profile'' of the gradient flow solution can change significantly depending on which ``flavor'' of inequality the cost satisfies. After a general theoretical analysis, we focus on fitting the CT-LQR policy optimization problem to the proposed framework, showing that, in fact, it can never satisfy a P\L I in its strongest form. We follow up our analysis with a brief discussion on the difference between continuous- and discrete-time LQR policy optimization, and end the paper with some intuition on the extension of this framework to optimization problems with $L_1$ regularization and solved through proximal gradient flows.

\end{abstract}

%%%%%%%%%%%%%%%%%%%%%%%%%%%%%%%%%%%%%%%%%%%%%%%%%%%%%%%%%%%%%%%%%%%%%%%%%%%%%%%%
\section{Introduction}

Recent advances in Artificial Intelligence (AI) and Machine Learning (ML) have rekindled interest in optimization theory, with many traditional results being revisited in light of the proposed techniques \cite{cui2024small,sontag_remarks_2022,de2024remarks,de2024convergence,mohammadi2021convergence,fatkhullin2021optimizing,bhandari_global_2022,eftekhari_training_2020}. In particular, the typical model-free formulation of many successful learning techniques motivates the study of gradient-based optimization methods, which are invaluable in understanding the training of neural networks and similar architectures, typically done through back-propagation algorithms. 

Gradient descent or, in continuous time, gradient flow, consists in searching for the argument $\x$ that minimizes the value of a given function $\cost[\x]$ by ``moving along'' the direction of steepest descent of the cost function. Theoretical guarantees are typically desirable, and in search of balancing generality and good properties, often in the optimization literature one deals with specific classes of optimization problems, such as convex optimization \cite{boyd2004convex} or linear programming \cite{chvatal1983linear}. In this paper, we will focus on optimization problems that satisfy (to different degrees) a Polyak-\L ojasiewisc inequality (P\L I), also known as the gradient dominance condition \cite{polyak1963gradient,karimi2016linear}. 

The P\L I is a staple in nonlinear optimization analysis, as, in its strongest form, it guarantees global exponential convergence of the gradient flow to the optimal solution of the problem \cite{karimi2016linear}. Furthermore, satisfying a P\L I globally (\gPLI) also guarantees strong robustness properties \cite{sontag_remarks_2022}. However, characterizing such a condition might not be possible for every optimization problem. In \cite{cui2024small} the authors noticed that more general conditions than the \gPLI can be proposed by using different classes of comparison functions so that different robustness results can be guaranteed. 

In particular, the problem of policy optimization for the linear quadratic regulator (LQR) motivates the discussion around weaker versions of the P\L I \cite{watanabe2025revisiting,fazel2018global,sun2021learning,hu2023toward,mohammadi_convergence_2022,fatkhullin2021optimizing,mohammadi2021lack,cui2024small,sontag_remarks_2022}. For the discrete-time version of the problem, in \cite{fazel2018global,sun2021learning,hu2023toward} the authors show that it satisfies a \gPLI, guaranteeing exponential convergence to the optimal feedback law for initialization in the stabilizing set of feedback matrices. However, so far in the literature for the continuous-time LQR policy optimization problem there is no characterization of a \gPLI \cite{mohammadi_convergence_2022,fatkhullin2021optimizing}, with the analysis in \cite{fatkhullin2021optimizing} indicating that, at least for the scalar case, the continuous-time LQR does not satisfy a \gPLI. 

In this work, we are interested in characterizing how generalizations of the P\L I affect the rate of convergence of the solution. We begin in Section \ref{sec:convopt} by revisiting common assumptions and their consequences regarding the convergence of the gradient flow. We then formally introduce the \gPLI and a few other weaker definitions, and discuss their differences and consequences to the convergence of the gradient flow solution. We next deepen the analysis by defining a new family of conditions closely related to, but more general than, the global P\L I. We discuss how these weaker conditions relate to each other and how they can characterize weaker forms of convergence than the \gPLI. Then, in Section \ref{sec:LQRPO}, we contextualize the theoretical analysis of this paper through the specific problem of the continuous-time LQR policy optimization problem. This problem has no guarantees of satisfying a \gPLI, and in fact we show that it can never satisfy such a condition. We characterize which sequences of points of the policy space result in an unbounded value for the gradient of the cost, and which result in a ``sub-exponential'' convergence profile for the solution. We follow up with a brief discussion on the difference between the continuous- and discrete-time LQR policy optimization, and finalize the paper in Section \ref{sec:conclusions} with a comment on possible links between the analysis of this paper and proximal gradient flow for optimization problems with $L_1$ regularizing terms.
All proofs are presented in the appendix for clarity.

%%%%%%%%%%%%%%%%%%%%%%%%%%%%%%%%%%%%%%%%%%%%%%%%%%%%%%%%%%%%%%%%%%%%%%%%%%%%%%%%
\section{Theoretical Setup}
\label{sec:convopt}

Along this paper, let $\re$, $\re_+$ and $\re_{++}$ denote the real, non-negative real, and strictly positive real numbers, respectively. Let $\PSD[n]$ and $\PD[n]$ be the set of positive semi-definite and positive definite $n$-by-$n$ matrices. 

A given function $\alpha:\re_+ \rightarrow\re$ is said to  be \emph{positive-definite} ($\mathcal{PD}$) if $\alpha(0)=0$ and $\alpha(x)>0$ for all $x\neq 0$. Similarly, $\alpha$ is said to be of class-$\mathcal{K}$ if it is continuous, positive-definite, and strictly increasing. Finally, $\alpha$ is of class-$\mathcal{K}_\infty$ if it is of class $\mathcal{K}$ and unbounded.

For a given function $\cost:\SSp\rightarrow\re$ bounded below, let $\minCost=\inf_{\x\in\SSp}f(x)$, and $\xs=\arg\inf_{\x\in\SSp}f(x)$.

\subsection{Optimization problems and gradient methods}

Let $\SSp$ be an open subset of an Euclidean space (the analysis in this paper can be generalized to manifolds, but we refrain from it for simplicity). Then, an optimization problem consists of searching for the value of an argument/parameter $\x\in\SSp$ that minimizes some cost function $\cost:\SSp\rightarrow\re$ (or minimizes the negative of a reward for maximization). Mathematically, we write such a problem as
\begin{equation}
    \label{eq:optprob}
    \begin{aligned}
        \underset{\x}{\textrm{minimize}} \quad & \cost[\x]\\
        \textrm{subject to} \quad & \x\in\SSp%\\
          %&\xi\geq0    \\
    \end{aligned}~~,
\end{equation}
which might have one, multiple, or no solution, requiring some assumptions about either the cost $\cost$ or the search space $\SSp$ to guarantee the existence and uniqueness of the solution. A common assumption is that of compactness of $\SSp$, and continuity of $\cost[\x]$ for $\x\in\SSp$, which would guarantee the existence of a minimum (and maximum) in $\SSp$ since $\cost[\SSp]:=\{\cost[\x]~|~\x\in\SSp\}$ would be compact. Usually, however, the search space $\SSp$ is not compact, requiring the adoption of an alternative set of assumptions, outlined next.
\begin{assumption}
    \label{asmp:bb&prp}
    The function $\cost$ is \emph{real analytic}, \emph{bounded below}, and \emph{proper} (\ie coercive).
\end{assumption}

It is easy to prove that Assumption \ref{asmp:bb&prp} guarantees the existence of a minimum $\minCost\in\cost[\SSp]$ attained at a set of points $\target:=\{\x\in\SSp ~|~ \cost[\x]=\minCost\}$. Notice that from the point of view of the optimization problem, any $\x\in\target$ is a valid solution of \eqref{eq:optprob}, as all result in the same value for the cost function. Then, the optimization problem can be thought of as finding any $x\in\target$. 

A natural candidate for solutions to the optimization problem is the set of critical points of $\cost$, i.e. the set of points $\CPts:=\{\x\in\SSp ~|~ \gradCost[\x]=0\}$. A common strategy for finding an $\x\in\CPts$ is ``moving the parameters along the direction of steepest descent of the function''. Mathematically and in continuous-time, this means imposing the following dynamics for the parameter $\x$
\begin{equation}
    \label{eq:gradflow}
    \dot \x = -\gradCost[\x],
\end{equation}
while in discrete time one would impose the following update law for $x_k$ for a small enough $h>0$
\begin{equation}
    \label{eq:graddesc}
    \x_{k+1} = \x_{k}-h\gradCost[\x_{k}].
\end{equation}

In this paper we focus on the continuous-time strategy, and one can easily verify that $\x$ is an equilibrium of \eqref{eq:gradflow} if and only if $\x\in\CPts$. Nonetheless, there is no a priori guarantee that a solution of \eqref{eq:gradflow} initialized in $\SSp$ will converge to a point in $\CPts$, much less in $\target$. So, we next look at what convergence guarantees Assumption \ref{asmp:bb&prp} allows us to derive, and what other assumptions can be made to improve such guarantees.

\subsection{Convergence guarantees and the Polyak-\L ojasiewicz inequality}
\label{ssc:ConvPLI}

Consider an optimization problem \eqref{eq:optprob} satisfying Assumption \ref{asmp:bb&prp}, then the following results hold:

\begin{lemma}
    \label{lem:precompact}
    For any \emph{proper} function $\cost:\SSp\rightarrow\re$, the solution of the gradient flow \eqref{eq:gradflow} initialized at any point $\x\in\SSp$ is \emph{precompact}.
\end{lemma}

\begin{theorem}[\L{}ojasiewicz's theorem \cite{lojasiewicz1984gradients}]
    \label{thm:Loj}
    Let $\cost:\SSp\rightarrow\re$ be a real analytic function and let $\sol:\re_+\times \SSp\rightarrow\SSp$, the solution of the initial value problem \eqref{eq:gradflow}, be precompact, i.e. $\sup_{t>0,\x_0\in\SSp}\|\sol[t,x_0]\|<\infty$. Then, for all $\x_0\in\SSp$, $\lim_{t\rightarrow\infty}\sol[t,x_0]\in\CPts$, that is, all solutions of the gradient flow \eqref{eq:gradflow} initialized in $\SSp$ converge to a critical point of the cost function $\cost$.
\end{theorem}

Lemma \ref{lem:precompact} and Theorem \ref{thm:Loj} guarantee that any solution of \eqref{eq:gradflow} initialized in $\SSp$ will converge to a critical point of the function $\cost$. However, $\x\in\CPts$ is only a necessary condition for the optimality of $\x$. In fact a point $\x\in\CPts$ can be either a local minimum, a local maximum, or a saddle-point of $\cost$. Regarding that, the following result can be stated (\cite{de2024remarks,de2024convergence} appendix A):% \LC{A reference is needed here.}

\begin{lemma}[\cite{de2024convergence,de2024remarks}]
    \label{lem:Conv2Min}
    Let $$\SSd:=\{\x\in\CPts~|~ \exists v \in \R^n ~\text{s.t.}~v^\top \nabla^2 \cost[\x]v<0\},$$ where $\nabla^2f$ is the Hessian of $\cost$, and let $\sol:\re_+\times \SSp\rightarrow\SSp$ be the solution of the initial value problem \eqref{eq:gradflow}. Then, the set of $x_0\in\SSp$ for which $\lim_{t\rightarrow\infty}\phi(t,x_0)\in\SSd$ has Lebesgue measure zero. In other words, the center-stable manifold of $\SSd$ has measure zero.
\end{lemma}

Notice that Lemma \ref{lem:Conv2Min} holds even if $\SSd$ is a continuous set, or the union of continuous sets. In fact, $\SSd$ need not be even compact, as long as the condition on the value of the Hessian holds for all of its elements. However, despite this result excluding any local maxima and ``strict" saddles from the result of a gradient flow solution to problem \eqref{eq:optprob} (with probability one), it is still not enough to guarantee the optimality of the gradient flow solution. Typically, other assumptions are added in the literature to ensure that $\lim_{t\rightarrow\infty}\sol[t,x_0]\in\target$, with, arguably, one of the most popular being the convexity of $\cost$. If $\cost$ is convex in $\SSp$, then $\CPts=\target$ and a gradient flow will eventually find the optimal solution. Furthermore, if $\cost$ is strongly convex, then a solution of \eqref{eq:gradflow} converges to $\target$ exponentially.

In this paper, we first review a condition for $\cost$ that is weaker than strong convexity but still ensures that a solution of \eqref{eq:gradflow} converges to $\target$ exponentially.

\begin{definition}[$\mu$-global Polyak-\L ojasiewicz inequality]
    \label{def:PLI}
    Given a fixed $\mu>0$, a function $\cost$ satisfies a \emph{$\mu$-global Polyak-\L ojasiewicz inequality} ($\mu$-\gPLI) if
    \begin{equation}
        \label{eq:PLI}
        \|\gradCost[\x]\|\geq \alpha\left(\cost[\x]-\minCost\right).
    \end{equation}
    with $\alpha(r) = \sqrt{\mu r}$, for all $\x\in\SSp$. Furthermore, $\cost$ is \gPLI if it is $\mu$-\gPLI for some $\mu>0$.
\end{definition}

The property in Definition \ref{def:PLI} (often written in the form
$\|\gradCost[\x]\|^2\geq \mu\left(\cost[\x]-\minCost\right)$)
has been the object of much recent study, and is a natural generalization of convexity (see \cite{karimi2016linear} for a through analysis of the relationship between convexity and the \gPLI). An immediate consequence of this property is that all critical points of $\cost$ must solve \eqref{eq:optprob}, \ie if $\cost$ is \gPLI, then $\CPts=\target$, which guarantees the optimality of gradient flow solutions. Further usefulness of this property lies in the fact that it guarantees exponential convergence of the solution of a gradient flow, as we show next.

\begin{definition}[$\mu$-global exponential stability]
    \label{def:glbexpst}
    Given a fixed $\mu>0$, the gradient flow \eqref{eq:gradflow} of $\cost$ is $\mu$-globally exponential stable ($\mu$-GES) if
    \begin{equation}
        \label{eq:glexpconv}
        \cost(\sol[t,x_0])-\minCost\leq (\cost[x_0]-\minCost)\mbox{e}^{-\mu t}.
    \end{equation}
    
    Furthermore, the gradient flow of $\cost$ is $GES$ if it is $\mu-GES$ for some $\mu>0$.
\end{definition}

\begin{lemma}
    \label{lem:expconvPLI}
    The gradient flow \eqref{eq:gradflow} of $\cost$ is $\mu$-GES if and only if $\cost$ satisfies a $\mu$-global P\L I.
\end{lemma}

Despite the good convergence properties associated with having an exponential upper-bound, proving that a cost function is $\mu$-global P\L I is not always possible, and in some relevant examples in the literature, the following weaker condition is characterized instead.

\begin{definition}[Semi-global P\L I]
    \label{def:sgPLI}
    A function $\cost$ satisfies a \emph{semi-global Polyak-\L ojasiewicz inequality} (\sgPLI) if, for every $\epsilon>0$, there exists a $\mu_\epsilon>0$ such that 
    %Let $\SSp_\epsilon$ be a sublevelset of $\cost$, \ie let $\SSp_\e:=\{\x\in\SSp~|~\cost[\x]-\minCost\le\e\}$. Then there exists some $\mu_\e>0$ for which the following inequality holds
    %
    \begin{equation}
        \label{eq:sgPLI}
        \|\gradCost[\x]\|\geq \alpha_{\e}(\cost[\x]-\minCost),
    \end{equation}
    with $\alpha_\e(r)=\sqrt{\mu_\e r}$, for all $x\in\SSp_\e:=\{\x\in\SSp~|~\cost[\x]-\minCost\le\e\}$.
\end{definition}

Similarly to satisfying a global P\L I, satisfying a semi-global P\L I guarantees that all critical points of $\cost$ must be global minima of \eqref{eq:optprob}, \ie $\CPts=\target$. In fact, at first glance global and semi-global P\L Is look very similar, with the latter also guaranteeing for any  $\e>0$, an exponential rate of convergence $\mu_\e$ for all initializations in the sublevel set $\SSp_\e$. However, their distinction becomes important when analyzing the rate of convergence of gradient flow solutions, as the following lemma illustrates.

\begin{lemma}
    \label{lem:PLIvssgPLI}
    If a function $\cost$ satisfies a \gPLI, then it also satisfies a \sgPLI. Alternatively, if a function $\cost$ satisfies a semi-global P\L I but not a global P\L I for any $\mu>0$, then there must exist a sequence $\{\e_i\}$, $i=1,2,\dots$ with $\e_i>0$ such that
    \begin{equation}
        \lim_{i\rightarrow\infty}\overline\mu_{\e_i}=0,
    \end{equation}
    where $\overline\mu_\e$ is the largest $\mu_\e$ that satisfies \eqref{eq:sgPLI} for a given $\e>0$.
\end{lemma}

Lemma \ref{lem:PLIvssgPLI} makes the distinction between satisfying a global and a semi-global P\L I clear: if the inequality is only semi-global, then there must be some ``unbounded sequence'' of points $\{x_i\}$, $x_i\in\SSp$, with $\lim_{i\rightarrow\infty}f(x_i)-\minCost=\infty$, for which the exponential rate of convergence $\overline\mu_{\cost[x_i]-\minCost}$ goes to zero as $i$ goes to infinity. Although intuitively this tells us that the convergence is not purely exponential globally, it is still unclear what the solution of the gradient-flow looks like. 

Furthermore, an even weaker version of the P\L I can be characterized as follows:

\begin{definition}[$\e$-local P\L I]
    \label{def:lPLI}
 Given a fixed $\e>0$, a function $\cost$ satisfies a \emph{$\e$-local Polyak-\L ojasiewicz inequality} ($\e$-{\lPLI}) if there exists some $\mu>0$ such that
    \begin{equation}
        \label{eq:lPLI-def}
        \|\gradCost[\x]\|\geq \alpha(\cost[\x]-\minCost),
    \end{equation}
    with $\alpha(r) = \sqrt{\mu r}$, for every $x\in\SSp_\e:=\{x\in\SSp~|~\cost[\x]-\minCost\leq\e\}$. Furthermore, $\cost$ is \lPLI if it is $\e$-\lPLI for some $\e>0$.
\end{definition}

Differently from global and semi-global P\L Is, a local P\L I only gives guarantees for fixed a neighborhood $\SSp_\e$ around the optimal value. Although inside such neighborhood the same guarantees are obtainable (exponential convergence, optimality of critical points, etc.), no guarantee exists outside of it in general.

Definitions \ref{def:PLI}, \ref{def:sgPLI} and \ref{def:lPLI} provide good granularity when analyzing the behavior of the gradient flow for different classes of cost functions, however further precision can be attained with the help of comparison functions, as we discuss next.

\subsection{A generalization of the P\L I }

The classic P\L I condition is originally formulated as in Definition \ref{def:PLI}, using the square root comparison function. However, practical cost functions often fail to satisfy this condition globally — an example being the  continuous-time LQR cost, which will be discussed later. 
This limitation motivated us in \cite{cui2024small} to propose a nonlinear version of the P\L I. By simply generalizing $\alpha(r)$ in \eqref{eq:PLI} from $\alpha=\sqrt{\mu r}$ to a positive definite function $\alpha \in \mathcal{PD}$, the convergence of the gradient flow can still be ensured. As an additional benefit, when the gradient flow in \eqref{eq:gradflow} is subject to additive noise, the error $f(\sol[t,\x_0]) - \minCost$ is bounded by an energy-like measure of the noise. This property is formally known as integral input-to-state stability (iISS) \cite{angeli2000characterization}. Furthermore, if $\alpha$ is strengthened to a class-$\mathcal{K}$ function, the error $f(\sol[t,\x_0]) - \minCost$ not only converges to zero in the absence of noise but also remains bounded when the noise is below a certain threshold, a property referred to as small-input ISS (siISS) \cite{cui2024small}. If $\alpha$ is further strengthened to be a class-$\mathcal{K}_\infty$ function, then the error $f(\sol[t,\x_0]) - \minCost$ converges to zero in the noise-free case and remains bounded under any bounded noise, which corresponds to the classical input-to-state stability (ISS) \cite{sontag1989smooth}. Clearly, the global P\L I belongs to the class of $\mathcal{K}_\infty$ functions, while the semi-global P\L I belongs to the class of $\mathcal{PD}$ functions.

{
For our analysis, we also introduce a new class of functions, called class-$\mathcal{K}_\mathrm{SAT}$ functions, which can be represented as
\begin{align}\label{eq:Ksatdef}
    \alpha(r) = \sqrt{\frac{ar}{b+r}} \quad \forall r \ge 0,
\end{align}
where $a,b>0$ are constants. 
}

{With this stabilished, the following definition summarizes the ``zoo" of the generalized inequalities based on different classes of comparison functions.  
}
{
\begin{definition}%[$\mathcal{K}_{\infty}$, $\mathcal{K}$, and $\mathcal{PD}$-P\L Is]
    \label{def:zoo}
    A function $\cost$ satisfies a \emph{class-$\mathcal{K}_{\infty}$ lower bound} (resp. class $\mathcal{K}_{\mathrm{SAT}}$,  $\mathcal{K}$, or $\mathcal{PD}$) if
    \begin{equation}\label{eq:classKPLI}
        \|\gradCost[\x]\|\geq \alpha\left(\cost[\x]-\cost[\xs]\right),
    \end{equation}
    for all $\x\in\SSp$, with $\alpha$ being a function of class-$\mathcal{K}_{\infty}$ (resp. class-$\mathcal{K}_{\mathrm{SAT}}$,  $\mathcal{K}$, or $\mathcal{PD}$). 
\end{definition}
}

\begin{figure}[t!]
    \centering
    \begin{tikzcd}[arrows=Rightarrow]
        \mbox{\gPLI} \arrow[r]\arrow[d]& \mathcal{K}_{\mathrm{SAT}} \arrow[d]\arrow[r]& \mbox{\sgPLI}\arrow[d]\arrow[r] & \mbox{{\lPLI}}\\ \mathcal{K}_\infty \arrow[r] \arrow[d,Leftrightarrow]& \mathcal{K} \arrow[r] \arrow[d,Leftrightarrow]& \mathcal{PD}\arrow[d,Leftrightarrow] \\ \mbox{ISS} & \mbox{siISS} & \mbox{iISS}
    \end{tikzcd}
    \caption{Diagram of the hierarchy between types of comparison functions and their relationship with the different types of Polyak-\L ojasiewicz inequality. Notice in particular that while all functions that satisfy a P\L I also have a class $\mathcal{K}_\infty$ lower-bound (as presented in definition \ref{def:zoo}), the converse is not necessarily true. Similarly, satisfying a sgPLI implies there exists a $\mathcal{PD}$ lower-bound, but the converse is not true. Also notice that the newly introduced class $\mathcal{K}_{\mathrm{SAT}}$ lower bound lies in between \gPLI and \sgPLI, and provides a better convergence guarantee. Finally, the \lPLI stands isolated in the graph, but one should note that it sustains the same convergence and robustness properties than the \gPLI, so long as the disturbance is small enough to keep the solution in a neighborhood of the optimum.}
    \label{fig:zoo}
\end{figure}
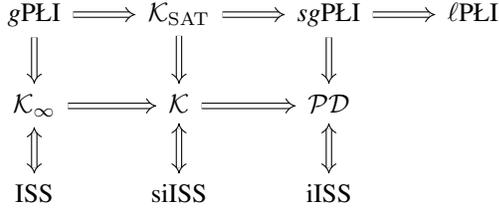

{Notice that there is a natural order between the comparison function and a natural relation between each of them and the previously defined different types of P\L I, all illustrated in Fig. \ref{fig:zoo}. In particular, notice that if the comparison function $\alpha$ is of class $\mathcal{K}_{\mathrm{SAT}}$, then it lies in between a \gPLI and a \sgPLI. %Interestingly if $\cost$ admits a class $\mathcal{K}_{\mathrm{SAT}}$ lower bound and has the gradient globally bound (\ie $\exists c>0$ such that $\|\gradCost[\x]\|\leq c$, $\forall x\in\SSp$), then the convergence rate is exponential within any sublevel set $\{\x~|~f(x)-f^*\le h\}$, but linear outside the sublevel set, as we formally show in the following results. }
Further relationships between the different classes can be established, as the following Lemma illustrates.

\begin{lemma}
    \label{lem:lpli-ksat}
    For a function $\cost$, the following two statements are equivalent
    \begin{enumerate}
        \item The function $\cost$ satisfies a class-$\mathcal{K}_\mathrm{SAT}$ lower bound.
        \item The function $\cost$ satisfies: a local P\L I; and a $\mathcal{PD}$ lower bound with comparison function $\alpha$ that is such that $$\liminf_{r\rightarrow\infty}\alpha(r)>0.$$
    \end{enumerate}
\end{lemma}

Satisfying different classes of lower bounds can also be shown to have consequences for the convergence rate of solutions. To better illustrate the effects of the different comparison functions, we next define a weaker form of convergence than the one from Definition \ref{def:glbexpst}.

\begin{definition}[global linear-exponential stability]
    \label{def:linexpconv}
    The gradient flow \eqref{eq:gradflow} of $\cost$ is globally linear-exponential stable (GLES) if there exists a $m>0$ such that for every $x_0\in\SSp$ and every $\underline t>0$ there exists a $\mu_{x_0,\underline t}>0$ for which
    \begin{align}
        \label{eq:linexp-bnd}
        (f(\phi(t,&\x_0)) - \minCost) \\&\le 
        \begin{cases} 
            (\cost(\sol[\underline t, \x_0])-\minCost)-m(t-\underline t) & \text{if } t\leq \underline t\\
            (\cost(\sol[\underline t, \x_0])-\minCost) e^{-\mu_{x_0,\underline t} (t - \underline{t})}  & \text{if } t > \underline{t}
        \end{cases}.\nonumber
    \end{align}
\end{definition}

From this definition, we can derive rigorous conditions for the solution to be GLES as follows:
{
\begin{lemma}
    \label{lem:linexp-bnd}
    The solution of the gradient flow \eqref{eq:gradflow} of a function $\cost$ is GLES if 
    \begin{itemize}
        \item The gradient $\gradCost$ is globally bounded;
        \item The function $\cost$ satisfies a class-$\mathcal{K}_\mathrm{SAT}$ lower-bound.
    \end{itemize}
\end{lemma}
}

\begin{figure}[t]
    \centering
    \includegraphics[width=1\linewidth]{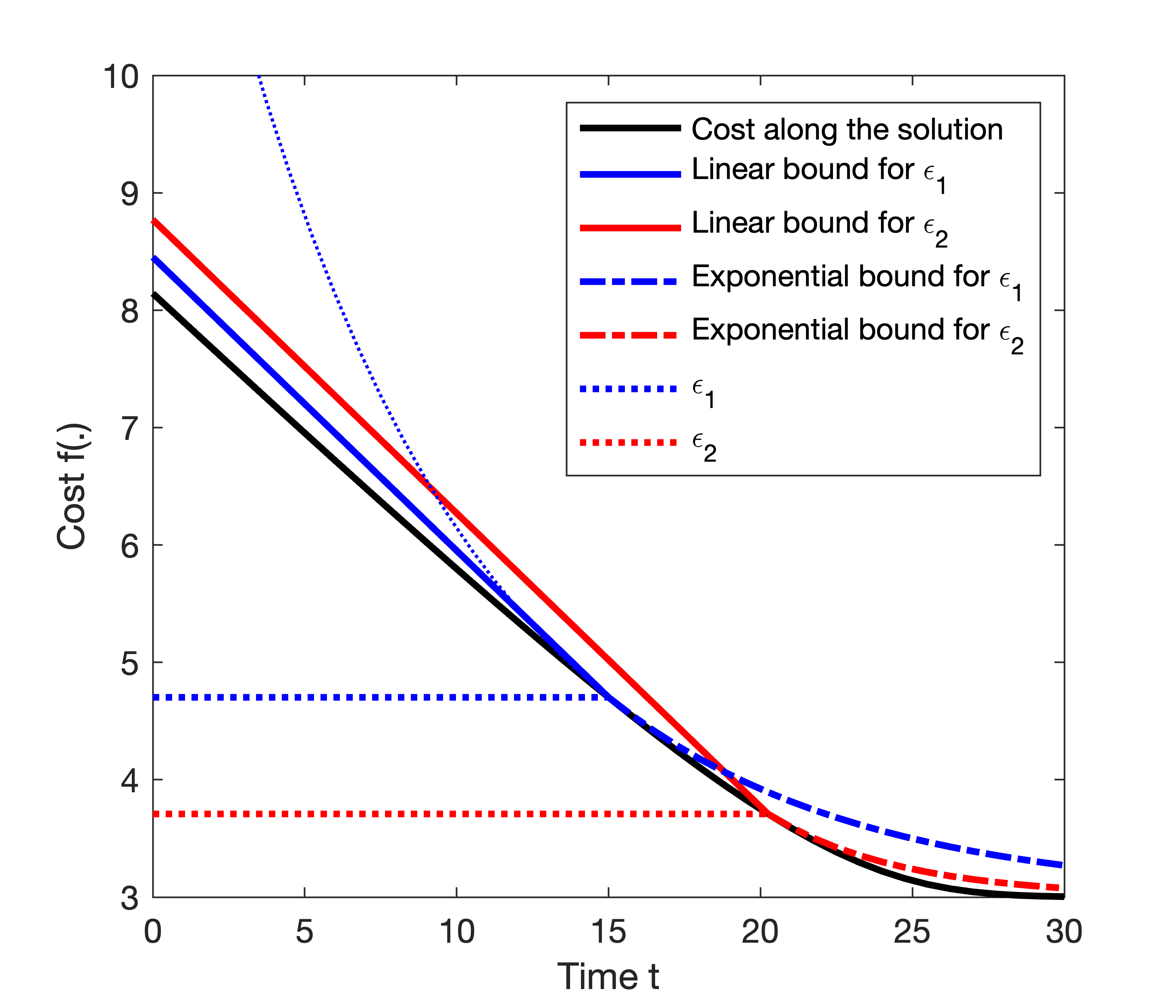}
    \caption{Gradient Flow trajectory for a cost $\cost$ that satisfies the conditions in Lemma \ref{lem:linexp-bnd}. Notice that for different given values of $\epsilon$, one can find a pair of line and exponential such that the line upper-bounds the solution while it is outside the levelset $\SSp_\e$, and the exponential upper-bounds it when it is inside the level-set. Furthermore, notice that smaller values of $\epsilon$ result in a tighter bound for the exponential section, while also on a looser bound for the linear section, while the opposite occurs when $\epsilon$ is larger. Finally, the thin dotted blue line illustrates how the exponential bound for $\e_1$ quickly becomes conservative if $t<\underline t$, pointing to linear-exponential as a tighter bound than purely exponential.}
    \label{fig:ThmIlustration}
\end{figure}

We illustrate the results from Lemma \ref{lem:linexp-bnd} in Fig. \ref{fig:ThmIlustration}. Observe that the true trajectory of the cost follows a ``linear-exponential'' profile, which means that for a given ``margin of error'' $\epsilon:=\cost(\sol[\underline t,x_0])-\underline f$, the solution is upper bounded by a line for values of the cost higher than $\epsilon$ and by and exponential for values smaller than $\epsilon$. Also notice from the figure that a smaller margin of error $\epsilon$ results in a worse upper bound for the linear part of the solution, but a tighter upper bound for the exponential part.

Next, we informally point out that for a given $x_0$ and $\underline t$, the exponential $(\cost(\sol[\underline t,x_0])-\minCost)\mbox{e}^{-\mu_{x_0,\underline t}(t-\underline t)}$ upper-bounds the solution $\cost(\sol[t,x_0])-\minCost$ for all time $t>0$, not only $t>\underline t$, however, for $t<\underline t$, the distance between the exponential upper-bound and the actual solution quickly grows. This can be noticed from Fig. \ref{fig:ThmIlustration} by looking at the thin dotted extension of the blue exponential for $t<\underline t = 15$. This means that the solution of this gradient flow is, in a sense, globally upper-bounded by an exponential. However, this bound is not tight for points much larger than $\phi(\underline t,x_0)$ (in the sense of $f(\phi(\underline t,x_0))$).

Moreover, notice that the conditions in Lemma \ref{lem:linexp-bnd} are only sufficient but not necessary. That is because if $\cost$ has a globally bounded gradient and satisfies both a local P\L I and a $\mathcal{PD}$ lower-bound, but is such that $\liminf_{r\rightarrow\infty}\alpha(r)=0$ for all positive-definite $\alpha$ for which $\cost$ satisfies \eqref{eq:classKPLI}, then there is still a linear-exponential upper-bound as described. However, in that case, the convergence of the solution can also be upper-bounded by a ``log-exponential'' function constructed similarly to how the linear-exponential bound is built in \eqref{eq:linexp-bnd}, and as a consequence, any linear function bounding the solution for $t<\underline t$ will also not be tight for points much larger than $\phi(\underline t,x_0)$ (in the sense of $f(\phi(\underline t,x_0))$).

Finally, notice a gap between GES and GLES: if $\cost$ does not satisfy a global P\L I, but also does not have a globally bounded gradient, then its solution is neither GES nor GLES. This is an important observation for the next section of the paper, as we will show that this is precisely the case for the LQR cost function.

In this section we provided tools for analyzing cost functions $\cost$ that satisfy Assumption \ref{asmp:bb&prp} and have one of the properties characterized in Definitions \ref{def:PLI}, \ref{def:sgPLI}, \ref{def:lPLI}, and \ref{def:zoo}. In the next section, we use these tools to analyze a very important example from the literature: the policy optimization problem for the Linear Quadratic Regulator (LQR) .

%%%%%%%%%%%%%%%%%%%%%%%%%%%%%%%%%%%%%%%%%%%%%%%%%%%%%%%%%%%%%%%%%%%%%%%%%%%%%%%%
\section{Applications to LQR policy optimization}
\label{sec:LQRPO}

Consider the following continuoue-time linear system:
\begin{align}
    \dot{x} &= Ax + Bu \label{eq:ltisys}
\end{align}
where $A\in\mathbb{R}^{n\times n}$, and $B\in\mathbb{R}^{n\times m}$ are the system matrices and are such that $(A,B)$ is controllable. Let $\Kstbl:=\{K\in\re[m\times n] ~|~ A-BK \text{ is Hurwitz}\}$. The objective is to determine an output feedback $u=-Ky$ with $K\in\Kstbl$ that minimizes 
\begin{align}
    \label{eq:costCTLQRbar}
    \bar J(K) = \mathbb{E}_{x_0\sim\mathcal{X}_0}\bigg[\int_0^\infty &x(t)^\top Q x(t) \\&+ x(t)^\top K^\top R K x(t) \, dt\bigg],\nonumber
\end{align}
with given matrices $R\in\PD[m\times m]$ and $Q\in\PD[n\times n]$, and for $x_0$ sampled from a probability distribution $\mathcal{X}_0$. It is well known that for linear systems minimizing $\bar J(K)$ is equivalent to minimizing the following cost function
\begin{equation}
    \label{eq:costCTLQR}
    J(K) = \trace{P_K},
\end{equation}
where
\begin{equation}
    \label{eq:LyapPCT}
    P_K(A-BK) + (A-BK)^\top P_K + K^\top RK + Q = 0,
\end{equation}
in the sense that a $K^*\in\Kstbl$ minimizes $J(\cdot)$ if and only if it minimizes $\bar J(\cdot)$. If matrices $A$ and $B$ are known, and the solution is initialized at a point $x_0$ sampled such that $\mathbb{E}(x_0x_0^\top) = \Sigma_0\in\PD[n]$, then one can find the optimal feedback matrix $K^*=R^{-1}B^\top P$ where $P$ solves the following Riccati equation
\begin{equation}
    \label{eq:RiccatiLQRCT}
    A^\top P+PA+PBR^{-1}B^\top P+Q=0.
\end{equation}
However, a popular formulation arises when one consider the case where the system matrices are not available, i.e. a model-free or policy optimization approach \cite{fazel2018global,hu2023toward,mohammadi_convergence_2022,watanabe2025revisiting,fatkhullin2021optimizing,cui2024small,de2024convergence,de2024remarks}. In that case, the optimal feedback matrix $K^*$ is obtained by following the negative direction of the gradient $\dot{K}=-\nabla J(K)$, with the gradient $\nabla J$ being given by \cite{rautert_computational_1997}:
\begin{equation}
    \label{eq:gradJCTLQR}
    \nabla J(K) = -2(B^\top P_K - RK) Y_K,
\end{equation}
%\LC{Should it be: $2(RK - B^\top P_K) Y_K$? }
%
where, for any $K\in\Kstbl$, $P_K$ is the solution of \eqref{eq:LyapPCT}, and $Y_K$ is the unique positive definite solutions of the following Lyapunov equation
\begin{equation}
    Y_K(A-BK)^\top + (A-BK)Y_K + I = 0. \label{eq:LyapYCT}
\end{equation}

In previous works in the literature \cite{mohammadi2021convergence,fatkhullin2021optimizing}, it was established that the solution for a gradient flow dynamics for solving this problem initialized inside $\Kstbl_a:=\{K\in\mathbb{R}^{m\times n}~|~J(K)\leq a\}$, satisfies a semi-global P\L I (Lemma 1 of \cite{mohammadi2021convergence}, and Theorem 3.16 of \cite{fatkhullin2021optimizing}), and in \cite{cui2024small} it was shown that it actually satisfies a class-$\Ksat$ lower-bound. Although, a priory, this does not mean that $J$ defined in \eqref{eq:costCTLQR} does not also satisfy a \gPLI, in the following section we will prove that this is actually the case, \ie $J(\cdot)$ can never satisfy a \gPLI.

\subsection{The LQR cost lies in the gap}

As mentioned, in this section we will show that the LQR cost $J(\cdot)$ has an unbounded gradient (thus not satisfying the conditions for Lemma \ref{lem:linexp-bnd}) while also provably not satisfying a \gPLI (thus not admitting a global exponential rate of convergence). 

We begin by formally defining ``high gain trajectories'' in the space of stable feedback matrices $\Kstbl$ and showing that, along any such high gain curve, the gradient is bounded, and thus the LQR cost can never satisfy a global P\L I. Then we show that for any sequence of matrices ``approaching the border of instability'' (in a sense to be formally defined) the gradient goes to infinity, proving that one cannot upper-bound the gradient globally in $\Kstbl$. 

Begin by defining ` gain curve'' in $\Kstbl$ as follows:

\begin{definition}
    \label{def:hgd}
    A matrix-valued function $\Ke:\re_+\rightarrow\Kstbl$, is called a \emph{high gain curve} of $\Kstbl$ if there exists an $\epsilon>0$ such that the eigenvalues of the closed-loop matrix are strictly in the left half-plane (LHP), \ie for all $\rho>0$, $$\max_i(\mathfrak{Re}(\lambda_i(A-B(\Ke(\rho)+K^*))))<-\epsilon,$$
    and the limit value of the cost is unbounded, \ie $\lim_{\rho\rightarrow\infty}J(K^*+\Kr)=\infty.$
\end{definition}

Furthermore, we can guarantee that any $(A,B)$ controllable has a high gain curve, as we show next:
\begin{lemma}
    \label{lem:hgd-exists}
    Let $A\in\re[n\times n]$ and $B\in\re[n\times m]$ be the system matrices for \eqref{eq:ltisys}, with $(A,B)$ being controllable. Then, there always exists a $\Ke:\re_+\rightarrow\Kstbl$ that is a high gain curve of $\Kstbl:=\{K\in\re[m\times n]~|~(A-BK)\text{ is Hurwitz}\}$, with the additional property that $B\Ke(\rho)$ is diagonalizable for all $\rho>0$. 
\end{lemma}

With this, we can state the following lemma regarding the behavior of the gradient along high gain trajectories:

\begin{lemma} \label{lem:BndGradJ}
    Let $\Ke:\re_+\rightarrow\Kstbl$ be a high gain curve of $\Kstbl$ and let there be a $\overline\rho>0$ such that $B\Kr$ is diagonalizable for for all $\rho>\overline\rho$. Then the limit
    \begin{equation}
        \lim_{\rho\rightarrow\infty}\nabla J(\Kr+K^*)\in\mathbb{R}^{m\times n},
    \end{equation}
    exists and is finite, \ie the norm of the gradient $\|\nabla J(\Kr+K^*)\|$ converges to a constant matrix as $J(\Kr+K^*)-J(K^*)$ goes to infinity.
\end{lemma}

Informally, Lemma \ref{lem:BndGradJ} proves the boundedness of the norm of the gradient along any high gain curve. As a consequence, it can never admit a class $\mathcal{K}_\infty$ P\L I lower-bound, as we formalize in the following corollary.

\begin{corollary}
    \label{cor:LQRnoPLI}
    There is no $\mu>0$ such that for all $K\in\mathbb{K}$ it holds that
    \begin{equation*}
        \|\nabla J(K)\| \geq \sqrt{\mu (J(K)-J(K^*))},
    \end{equation*}
    \ie, the cost function $J$ defined in \eqref{eq:costCTLQR} can never satisfy a \gPLI. 
\end{corollary}

% \begin{remark}
%     Notice that the result from Corollary \ref{cor:LQRnoPLI} is specially surprising if one considers that results from the literature characterize Assumption \ref{def:PLI} for the model-free discrete-time LQR \cite{fazel}. This distinction between the two formulations for the LQR raises questions regarding the role of sampling in the convergence of the gradient flow, but are not investigated in this paper due to limited space.
% \end{remark}

From these results, one would think that $J$ has a globally bounded gradient norm, since Lemma \ref{lem:BndGradJ} shows that the gradient converges to a constant matrix for any high gain curve. However, let $\bdKstbl := \{K\in\re[m\times n] ~|~ \lambda_{\max}(A-BK)=0\}$ be the border of stability, \ie the values of $K$ for which {$(A-BK)$ has at least one eigenvalue on the imaginary axis}, and notice that the value of the gradient \eqref{eq:gradJCTLQR} explodes to infinity as $K$ approaches $\bdKstbl$, as we show in the following lemma.

\begin{lemma}
    \label{lem:GradJUnbd}
    For any $\overline{K}\in\bdKstbl$, let $\{K_i\}$ for $i=1,2,\dots$ be a sequence of matrices $K_i\in\Kstbl$ such that $\lim_{i\rightarrow\infty}K_i=\overline K$, then $$\lim_{i\rightarrow\infty}\|\nabla J(K_i)\|_F=\infty.$$
\end{lemma}

As a consequence of this fact, gradient $\nabla J(\cdot)$ does not admit a global upper bound, and thus does not satisfy the conditions for Lemma \ref{lem:linexp-bnd}.

The fact that the continuous-time LQR cost neither satisfies a \gPLI, nor has a globally bounded gradient, makes it hard to provide tight global convergence rate estimates to the policy-optimization algorithm. In fact, the solution can be either exponential or linear-exponential, depending on which region of the state-space it is initialized at. Fortunately, since the gradient is only unbounded for ``bounded directions'' in the border of the space $\Kstbl$, we can enforce boundedness of the norm of the gradient on a subset of $\Kstbl$ as follows:

\begin{lemma}
    \label{lem:LQRClassKPLIub}
    For a $\delta>0$, let $\Ksstbl[\delta] :=\{K\in\Kstbl ~|~ A-BK+\delta I \text{ is Hurwitz}\}$, then there exists $g>0$ such that $\|\nabla J(K)\|<g$ for all $K\in\Ksstbl[\delta]$.
\end{lemma}

The results presented so far illustrate, through the LQR policy optimization problem, the value of understanding exactly what kind of ``P\L I-like'' condition the cost function in question satisfies. To complement our analysis so far, and in hopes of illustrating the different possible behaviors for the solution of the policy optimization for the LQR, we next provide an analysis of the single-input single-state/output LQR case.

\subsection{Convergence analysis of the scalar LQR policy optimization}

For the scalar case, the continuous-time system dynamics is given by \eqref{eq:ltisys} with $A=a\in\re$, and $B=b=1$, the later being assumed without loss of generality, since for different values of $b$, its magnitude could simply be included in the magnitude of the considered input signal $u$. The weighting matrices of the LQR cost are $Q=q \in \R$ and $R=r\in \R$.

Then, the cost and its gradient can be computed for the scalar case as
\begin{equation*}
    J(k)=p,
\end{equation*}
and
\begin{equation*}
    \partial J(k) = -2(p-rk)\ell,
\end{equation*}
where $p$ and $\ell$ solve \eqref{eq:LyapPCT} and \eqref{eq:LyapYCT}, respectively, for scalar parameters. From this, we can recover the results of Lemma \ref{lem:BndGradJ} for the scalar case. Notice that
\begin{align*}
    \partial J(k) &= -\frac{2rk(a-k)+(rk^2+q)}{2(a-k)^2}\\
    &= -\frac{2r(\frac{a}{k}-1)+(r+\frac{q}{k^2})}{2(\frac{a}{k}-1)^2},
\end{align*}
which implies that $\lim_{k\rightarrow +\infty}\partial J(k)= r/2$, indicating linear convergence as the magnitude of the feedback gain $k$ increases while still keeping $a-k<0$. This linear convergence, however, becomes less noticeable as the solution approaches the optimum feedback value $k^*$ and instead exponential convergence is observed. To show that we compute 
\begin{equation*}
    m(k) = \frac{\|\partial J(k)\|^2}{J(k)-J(k^*)}
\end{equation*}
for $k=k^*+\epsilon$, and verify that as $\epsilon\rightarrow 0$, $m(k)\rightarrow c$ for some positive constant $c$. That is done in the following proposition

\begin{proposition}
    For the scalar LQR, let $k^*$ be the value of the feedback gain that minimizes the cost $J(k)$. Then we have that 
    \begin{align}
        &J(k^*+\epsilon)-J(k^*) = \ell r \epsilon^2=:\delta(\epsilon)\\
        &\partial J(k^*+\epsilon) = \ell({-}\delta(\epsilon)+r\epsilon)\\
        \label{eq:me}
        &m(k^*+\epsilon) = r\ell(\ell^2\epsilon^2 {-}2\ell\epsilon+1).
    \end{align}
    These computations allow us to conclude that $m(k^*) = r\ell^*>0$, characterizing exponential convergence near $k^*$. Furthermore, $\lim_{\epsilon\rightarrow a-k^*}m(k^*+\epsilon)=\infty$, indicating that the convergence rate explodes for values of $k$ in the boundary of stability.
\end{proposition}

The result above was already known from the general case, but becomes clearer when derived for the scalar case, independent of matrix equalities and different ``high gain trajectories''. Furthermore, the simpler form of the scalar case makes it easier to analyze numerical results.

\begin{figure}[t!]
    \centering
    \includegraphics[width=0.49\linewidth]{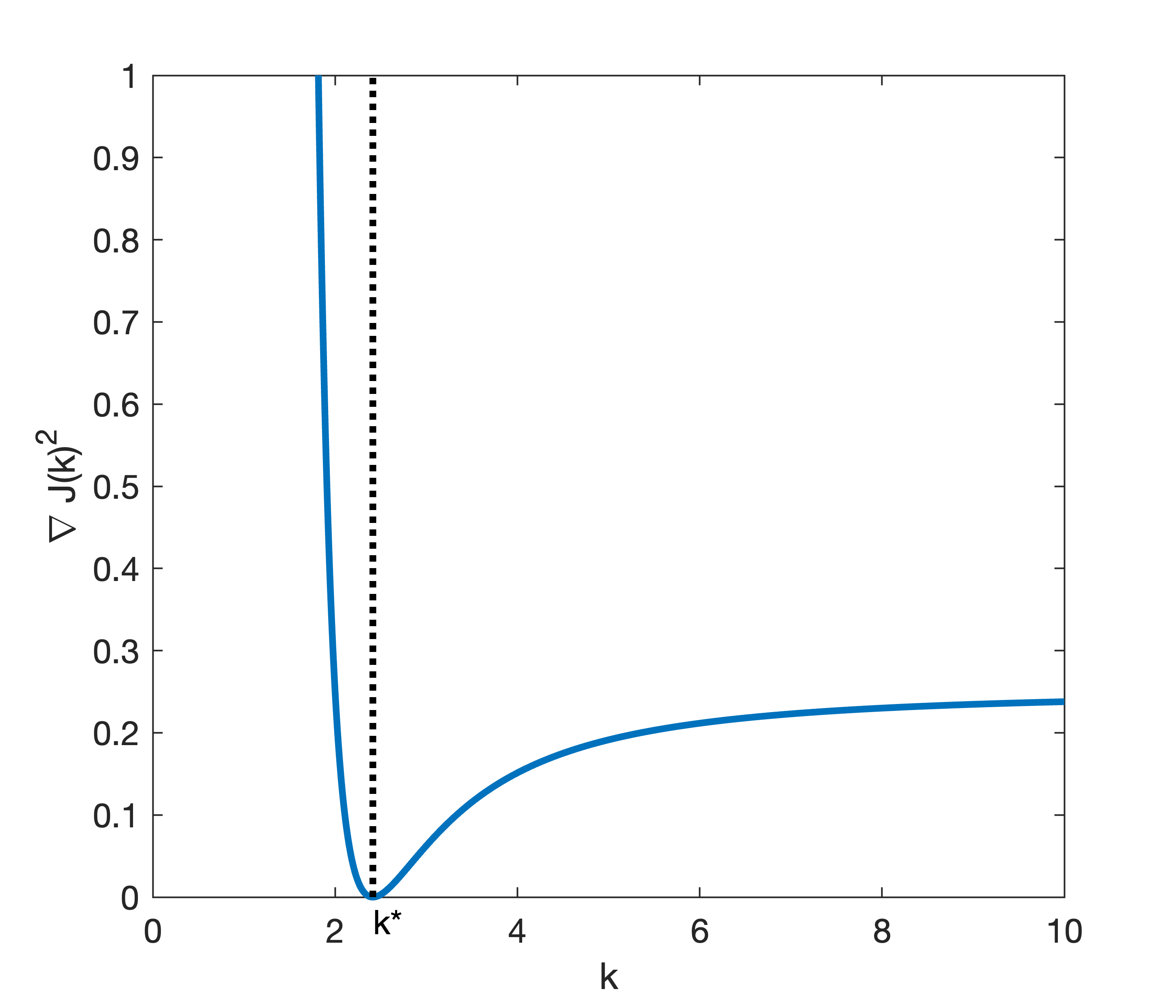} 
    \includegraphics[width=0.49\linewidth]{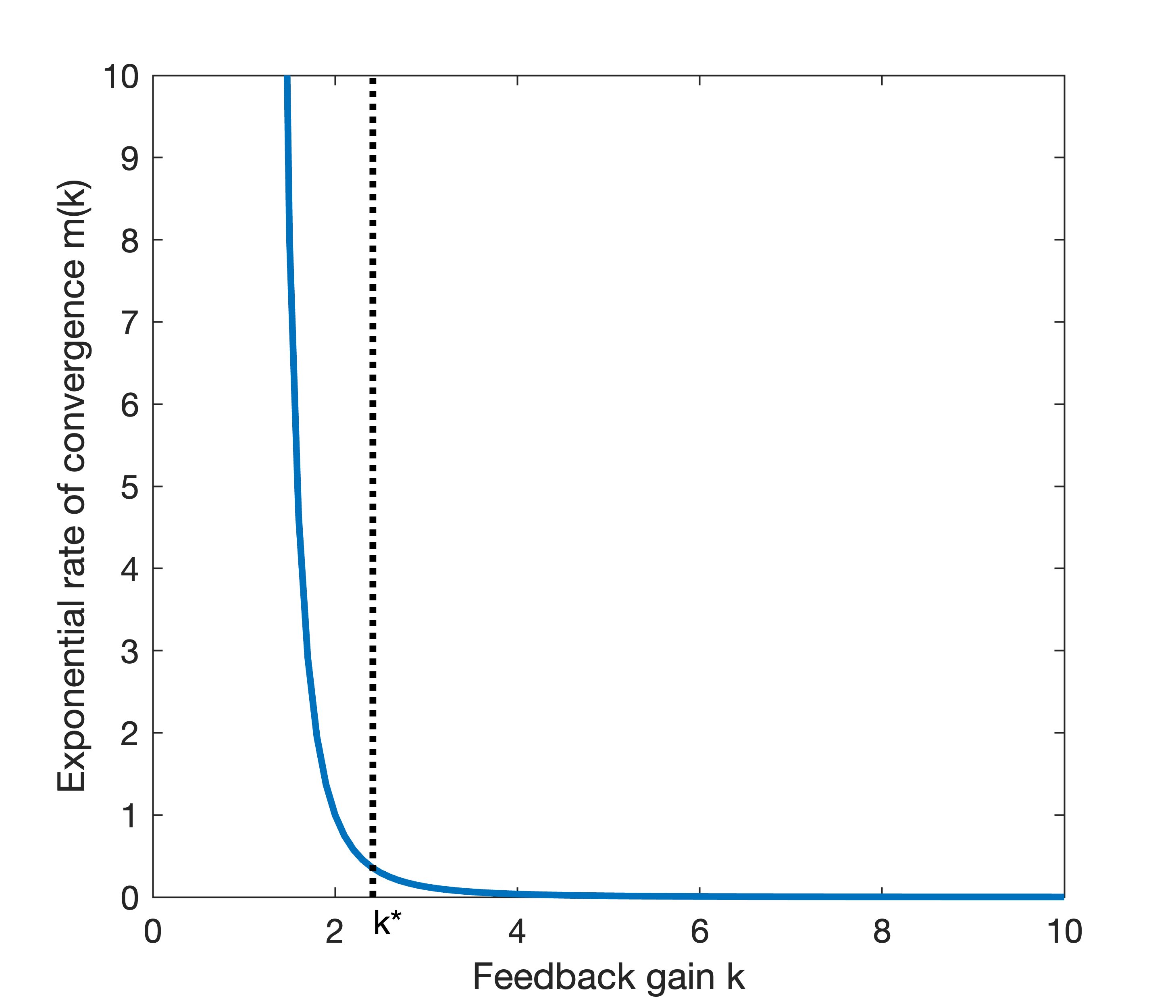} \\ (a)~~~~~~~~~~~~~~~~~~~~~~~~~~~~~~~~(b)
    \caption{Illustration of the dual behavior of the LQR cost function \eqref{eq:costCTLQR} fort he scalar case. In (a) we plot the squared norm of the gradient $\|\nabla J\|^2$ as a function of the feedback gain $k$, while in (b) we plot the largest exponential rate of convergence $m(k)$ as defined in \eqref{eq:me}, and in both plots the dotted vertical line indicates the optimal $k^*$. Notice that for $k>k^*$ (right side of the dotted line) the gradient is bounded above, and $m(k)$ quickly goes to zero, while for $k<k^*$ (left side of the dotted line) the value of the gradient, and the best exponential rate of convergence both quickly diverge to infinity as $k$ approaches the border of instability.}
    \label{fig:gradJscalar}
\end{figure}

Take $a=q=r=1$ and notice from Fig. \ref{fig:gradJscalar} that the gradient $\partial J(k)$ behaves completely differently if $k>k^*$ or if $k<k^*$. However, notice that if we restrict the domain from $[1,\infty]$ to $[1+\epsilon,\infty]$, the value of the gradient is now globally bounded. This is the intuition behind Lemma \ref{lem:LQRClassKPLIub}.

Furthermore, the linear-exponential convergence behavior described in Lemma \ref{lem:linexp-bnd} becomes very evident for the scalar LQR if $k$ is initialized larger than $k^*$, as can be seen in Fig. \ref{fig:gradflowsolscalarLQR}a. If, however, the solution is initialized near the border of instability, the convergence is much closer to a decreasing exponential, as evident in Fig. \ref{fig:gradflowsolscalarLQR}b.

\begin{figure}[t!]
    \centering
    \includegraphics[width=0.49\linewidth]{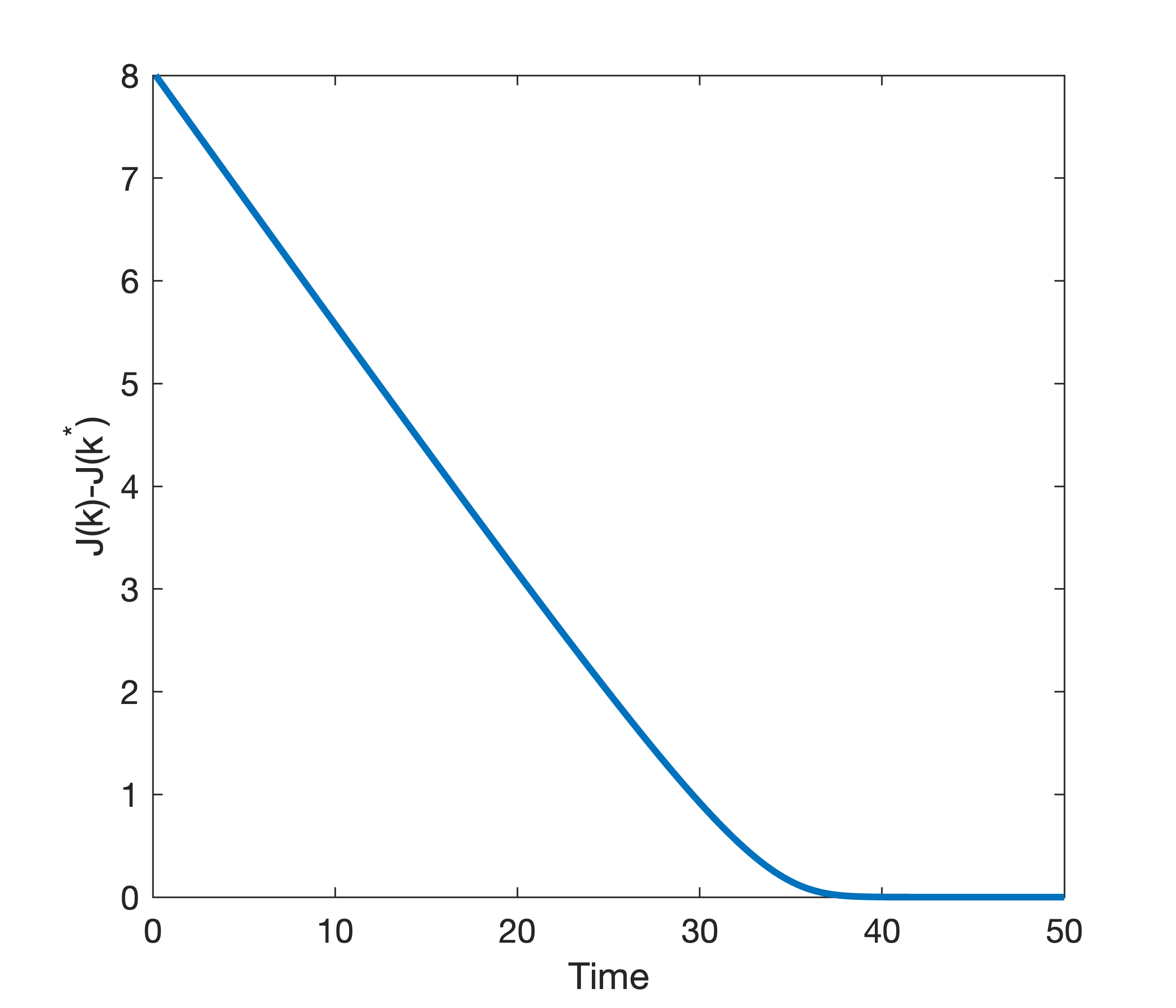} 
    \includegraphics[width=0.49\linewidth]{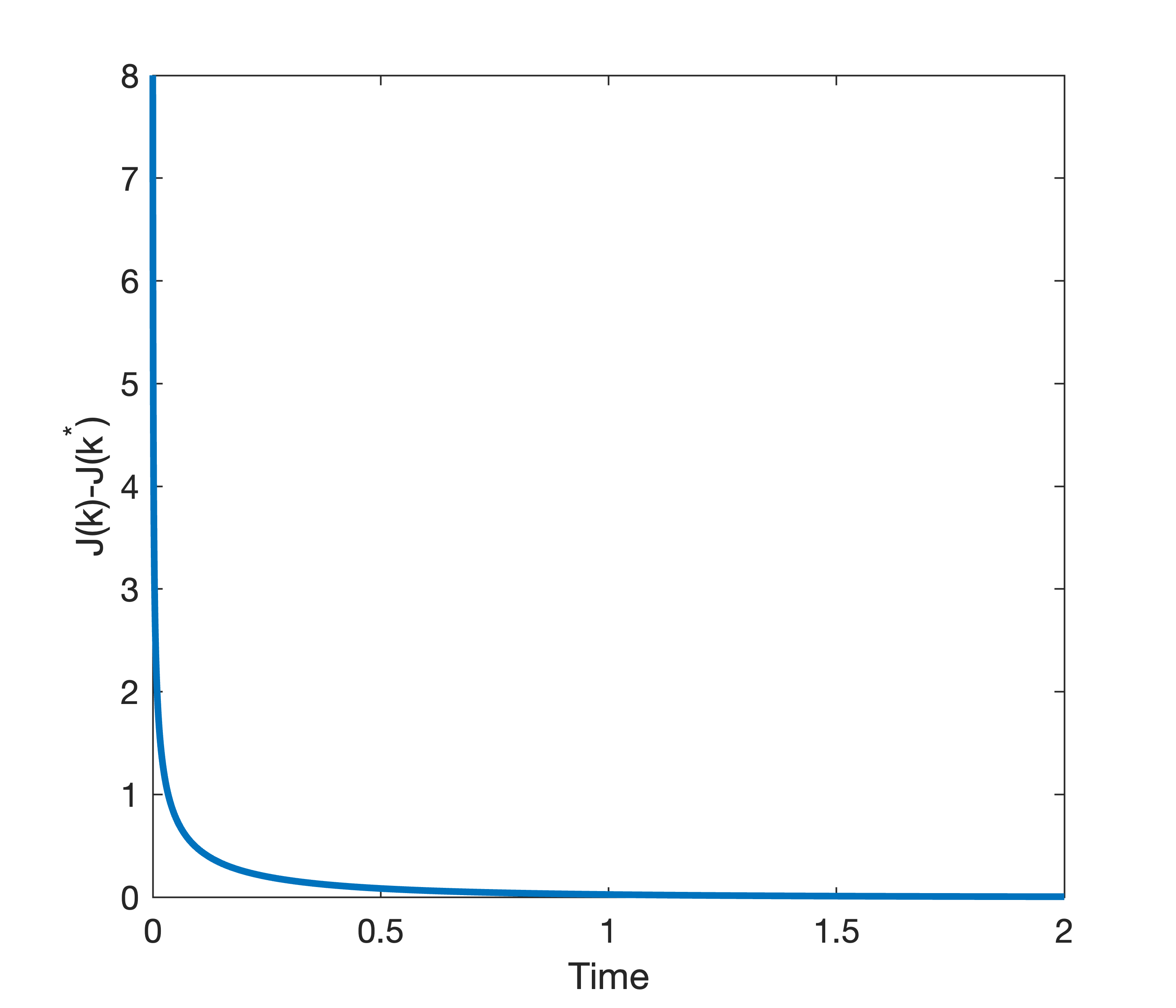} \\ (a)~~~~~~~~~~~~~~~~~~~~~~~~~~~~~~~~(b)
    \caption{Simulation results for the gradient flow of the scalar LQR policy optimization with $a=r=q=1$. Both simulations (a) and (b) were initialized such that $J(k(0))-J(k^*)\approx 8$, however (a) was initialized for $k(0)>k^*$ and (b) for $k(0)<k^*$. Notice that in (a) the convergence is ``linear exponential'' as described in Section \ref{ssc:ConvPLI} since, as can be seen in Fig. \ref{fig:gradJscalar}, for $k>k^*$, $\nabla J(k)$ is bounded above. In (b), on the other hand, the convergence is much quicker and exponential, due to the fact that for $k\in[a,k^*]$, the exponential rate of convergence $m(\epsilon)$ defined in \eqref{eq:me} is bounded away from zero. }
    \label{fig:gradflowsolscalarLQR}
\end{figure}

\subsection{Comments on the difference between continuous and discrete-time LQR policy optimization}

We conclude the analysis of this paper with a brief overview of the behavior of the discrete-time LQR policy optimization problem. This scenario is studied in different papers in the literature \cite{fazel2018global,sun2021learning,hu2023toward,mohammadi_convergence_2022} and for this problem, a global P\L I is characterized (see, for example, Lemma 1 in \cite{hu2023toward}). This is surprising since, by Corollary \ref{cor:LQRnoPLI}, the continuous-time LQR policy optimization can never admit a global P\L I.

Some intuition behind this difference can be obtained by looking at the Euler discretization of the scalar continuous case with step-size $h>0$, and with $R_d = hr$ and $Q_d=hq$. For this problem, we can define $m_d(k_d^*+\e)$ similarly to how it was done for the continuous-time case in \eqref{eq:me}. %We omit the explicit derivation here due to space limitations but have it in the extended version of this paper in arXiv \cite{arxiv}.
%
% \begin{align*}
%     m_d(k^*&+\e,h)  \\ =& (hp^*+r)\ell(\epsilon)\left(h^2+2\ell(\epsilon)(1+h(a+\kdst+\epsilon))\epsilon\\&+\epsilon^2\ell^2(1+h(a+\kdst+\epsilon))^2\right).
% \end{align*}
%

Furthermore, notice that the feedback gain $k_d$ is bounded between $a<\kdst+\epsilon<(2+ha)/h$, which implies that $\Kstbl$ is compact in the discrete-time. Moreover, upon explicit computation of $m_d(k_d^*+\e)$, one can check that for any $h>0$, $m_d(k_d,h)\rightarrow\infty$ if either $k_d\rightarrow a$ or $k_d\rightarrow (2+ha)/h$, which implies that $m_d(\cdot, h)$ must admit a minimum value $\underline{m_d}(h)$ attained at some point $\underline{\kd}\in\Kstbl$, \ie $\underline{m_d}(h) = m_d(\underline{\kd},h)$.

\begin{figure}[t!]
    \centering
    \includegraphics[width=0.49\linewidth]{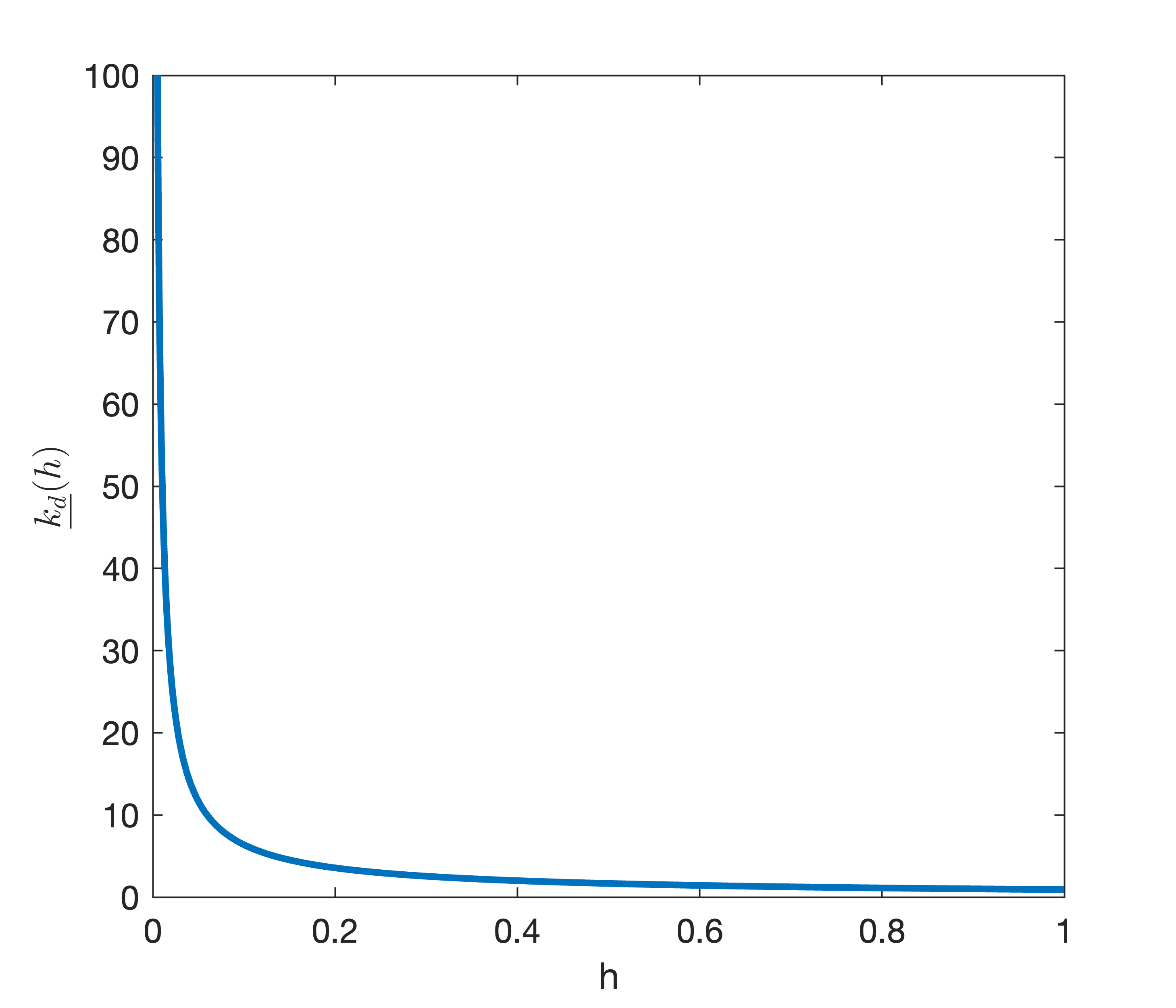} 
    \includegraphics[width=0.49\linewidth]{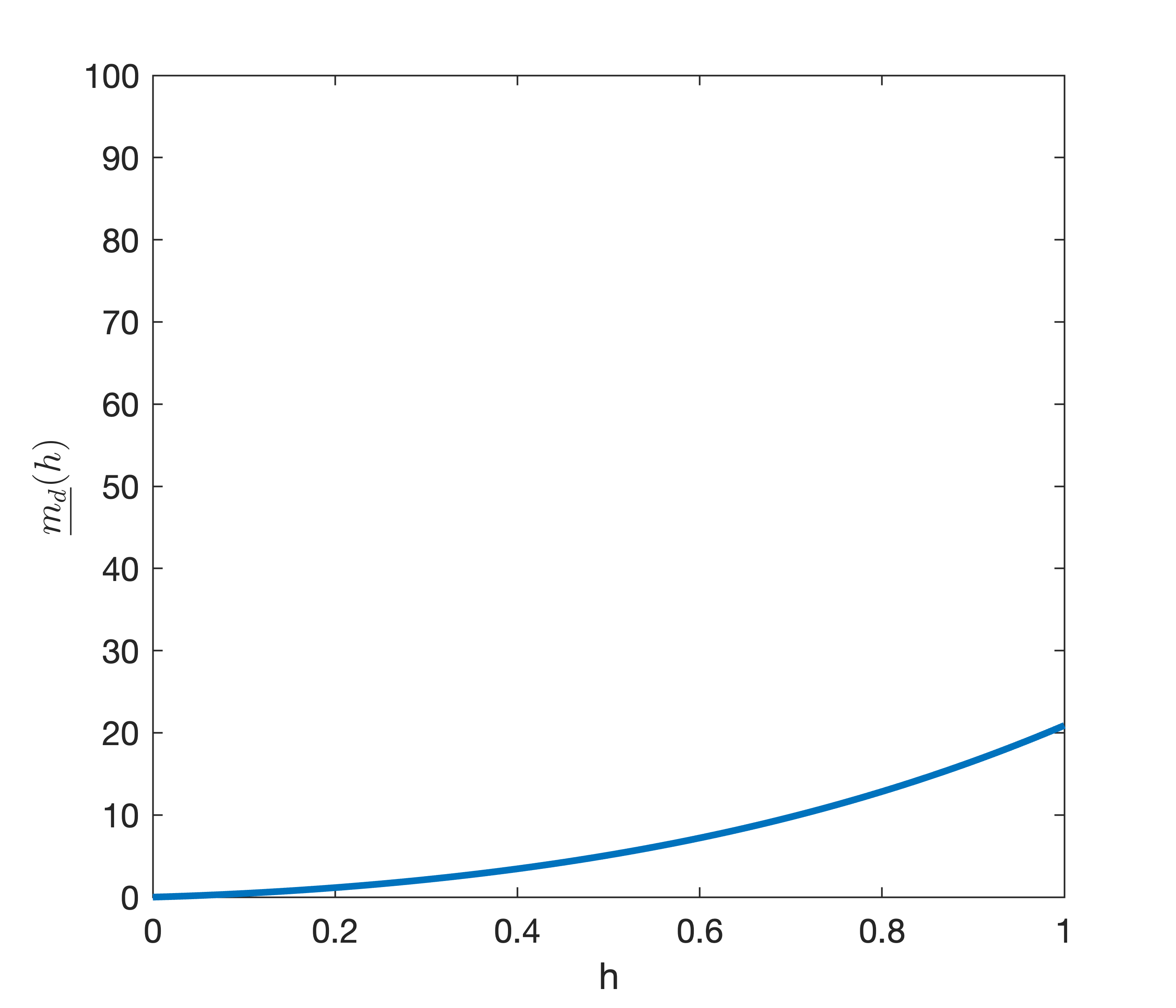} \\ (a)~~~~~~~~~~~~~~~~~~~~~~~~~~~~~~~~(b)
    \caption{Visualization of how the discretization step affects the global exponential rate of convergence for the scalar discrete-time LQR policy optimization problem. Notice that for any discretization step $h>0$, there exists a $\mu=\underline{m_d}(h)$ that provides global exponential convergence guarantees to the solution, however, that rate of convergence goes to zero as $h$ goes to zero, which is compatible with the observation that CT LQR does not have a global exponential rate of convergence.}
    \label{fig:dtlqr-disceffect}
\end{figure}

With these established, we pick $a=q=r=1$ and plot $\underline{\kd}(h)$ and $\underline{m_d}(h)$ in Fig. \ref{fig:dtlqr-disceffect}. Notice that as $h\rightarrow 0$ (\ie as we approach the CT LQR policy optimization problem), the value $\underline{k_d}(h)$ at which $m_d(\kd,h)$ is minimized goes to infinity, while $\underline{m_d}(h)$ goes to zero.

%%%%%%%%%%%%%%%%%%%%%%%%%%%%%%%%%%%%%%%%%%%%%%%%%%%%%%%%%%%%%%%%%%%%%%%%%%%%%%%%

% \section{A System with strictly PD PLI upper-bound(?)}

%%%%%%%%%%%%%%%%%%%%%%%%%%%%%%%%%%%%%%%%%%%%%%%%%%%%%%%%%%%%%%%%%%%%%%%%%%%%%%%%

\section{Conclusions and future works}
\label{sec:conclusions}

In this paper, we present a brief overview of convergence guarantees for gradient methods in optimization problems. We revisited the Polyak-\L ojasiewicz inequality (\ie gradient dominance condition) and observed how slight changes in its characterization can imply significant changes to the convergence of the gradient flow solution. This motivated the introduction of nonlinear comparison functions as a way of characterizing the behavior of the solution, which we supported with a result that gives conditions for the solution to present a ``linear-exponential'' behavior in Lemma \ref{lem:linexp-bnd}.

The paper follows up with a scenario where the traditional P\L I condition does not hold: the continuous-time model-free linear quadratic regulator problem. For this problem we showed that it presents neither globally exponential, nor globally linear-exponential convergence behavior, but a mixture of both depending on how close the solution is initialized to the border of instability. Despite that, we show in Lemma \ref{lem:BndGradJ} that for any ``high gain curve'' in the space of stabilizing feedback matrices, the gradient is upper-bounded, which allowed us, through Lemma \ref{lem:LQRClassKPLIub}, to characterize global linear-exponential convergence behavior through a judicious restriction of the optimization search space. We then illustrate our results through numerical simulations of the scalar case, where the two regions of the parameter space are clearly defined.

\subsection{A brief comment on lasso and proximal gradient}

To finish the paper, we offer a brief comment illustrating the consequences of adding an $L_1$ regularization term to the cost in \eqref{eq:optprob}, and solving it through proximal gradient flow instead of gradient flow (we refer the reader to \cite{hassan2021proximal,gokhale2024proximal} for a complete characterization of the expressions used here). Assume $\SSp\subseteq\re$ for simplicity (scalar parameter), then the optimization problem with an $L_1$ regularization is written as follows
\begin{equation}
    \label{eq:optproblasso}
    \begin{aligned}
        \underset{\x}{\textrm{minimize}} \quad & \cost[\x]+|\x|\\
        \textrm{subject to} \quad & \x\in\SSp%\\
          %&\xi\geq0    \\
    \end{aligned}~~,
\end{equation}
which cannot be solved through gradient flow due to the non-differentiability of $|\cdot|$. Instead, a solution is found through proximal gradient flow, which can be written as
%
% kdot = -(k-prox_{m*g}(k-m*Df(k))) = -k+sign(k-m*Df(k))*max(|k-m*Df(k)|-m,0)
\begin{equation*}
    \dot x = -(\x-\prox_{|\cdot |}(\x-\gradCost[\x]),
\end{equation*}
where $\prox_{|\cdot|}(x)$ is itself the solution of an optimization problem, but which has the following closed-form solution since we use $L_1$-norm:
\begin{equation*}
    \prox_{|\cdot|}(\x) = \sign(\x)\max(|\x|-1,0),
\end{equation*}
which after substituting to the proximal gradient flow results in
\begin{equation*}
    \dot x = -\x +\sign(x-\gradCost[\x])\max(|x-\gradCost[\x]|-1,0).
\end{equation*}

At this point we note informally that if $\cost$ satisfies the conditions in Lemma \ref{lem:linexp-bnd} (specifically global boundedness of the gradient), then there must exist a $\epsilon>0$ such that for all $\x$, if $\cost[\x]-\minCost\geq\epsilon$, then $\x-\gradCost[\x]\geq1$. Therefore, for initializations ``large enough'' (in the sense that $f(\x_0)-\minCost\geq \epsilon$) the proximal gradient flow simplifies to
\begin{equation*}
    \dot{\x} = -\gradCost[\x]-1,
\end{equation*}
Alternatively, if $\SSp$ is such that any sequence $\{\x_i\}$ approaching a finite point in its boundary is such that $\lim_{i\rightarrow\infty}\nabla f(x_i)=\infty$, then at a point $\x\in\SSp$ ``close enough'' to the boundary, $\gradCost[\x]\gg\x$ which simplifies the proximal gradient flow expression to
\begin{equation*}
    \dot{\x} = -\gradCost[\x]+1.
\end{equation*}

Neither case changes the global behavior predictions given in this paper for the solution while it remains far away from $\target$ (in some sense). In fact, one can argue that if either $|\x|\gg1$ or $|\gradCost[\x]|\gg1$ then $|x-\gradCost[\x]|-1<0$ only if $x\approx \gradCost[\x]$ (in the sense that $|\x-\gradCost[\x]|$ is much smaller than either $|\x|$ or $|\gradCost[\x]|$), which means that $\dot\x = -\x$ is an approximation of $\dot \x = -\gradCost[\x]$. However, notice that it is not straightforward to provide an equivalent local analysis. 

Beyond its practical relevance, this observation motivates technically future works on proximal gradient flows, and its equivalent P\L Is, which we believe are a natural follow-up to this publication.

%%%%%%%%%%%%%%%%%%%%%%%%%%%%%%%%%%%%%%%%%%%%%%%%%%%%%%%%%%%%%%%%%%%%%%%%%%%%%%%%

\bibliographystyle{ieeetran}
\begin{spacing}{0.86}
\bibliography{LQRPOConv_Bibliography}

\begin{thebibliography}{10}
\providecommand{\url}[1]{#1}
\csname url@rmstyle\endcsname
\providecommand{\newblock}{\relax}
\providecommand{\bibinfo}[2]{#2}
\providecommand\BIBentrySTDinterwordspacing{\spaceskip=0pt\relax}
\providecommand\BIBentryALTinterwordstretchfactor{4}
\providecommand\BIBentryALTinterwordspacing{\spaceskip=\fontdimen2\font plus
\BIBentryALTinterwordstretchfactor\fontdimen3\font minus \fontdimen4\font\relax}
\providecommand\BIBforeignlanguage[2]{{%
\expandafter\ifx\csname l@#1\endcsname\relax
\typeout{** WARNING: IEEEtran.bst: No hyphenation pattern has been}%
\typeout{** loaded for the language `#1'. Using the pattern for}%
\typeout{** the default language instead.}%
\else
\language=\csname l@#1\endcsname
\fi
#2}}

\bibitem{cui2024small}
L.~Cui, Z.-P. Jiang, and E.~D. Sontag, ``Small-disturbance input-to-state stability of perturbed gradient flows: {A}pplications to {LQR} problem,'' \emph{Systems \& Control Letters}, vol. 188, p. 105804, 2024.

\bibitem{sontag_remarks_2022}
\BIBentryALTinterwordspacing
E.~D. Sontag, ``\BIBforeignlanguage{en}{Remarks on input to state stability of perturbed gradient flows, motivated by model-free feedback control learning},'' \emph{\BIBforeignlanguage{en}{Systems \& Control Letters}}, vol. 161, p. 105138, Mar. 2022. [Online]. Available: \url{https://linkinghub.elsevier.com/retrieve/pii/S0167691122000056}
\BIBentrySTDinterwordspacing

\bibitem{de2024remarks}
A.~C.~B. De~Oliveira, M.~Siami, and E.~D. Sontag, ``Remarks on the gradient training of linear neural network based feedback for the {LQR} problem,'' in \emph{2024 IEEE 63rd Conference on Decision and Control (CDC)}, 2024, pp. 7846--7852.

\bibitem{de2024convergence}
A.~C.~B. de~Oliveira, M.~Siami, and E.~D. Sontag, ``Convergence analysis of overparametrized {LQR} formulations,'' \emph{arXiv preprint arXiv:2408.15456}, 2024.

\bibitem{mohammadi2021convergence}
H.~Mohammadi, A.~Zare, M.~Soltanolkotabi, and M.~R. Jovanovi{\'c}, ``Convergence and sample complexity of gradient methods for the model-free linear-quadratic regulator problem,'' \emph{IEEE Transactions on Automatic Control}, vol.~67, no.~5, pp. 2435--2450, 2021.

\bibitem{fatkhullin2021optimizing}
I.~Fatkhullin and B.~Polyak, ``Optimizing static linear feedback: {G}radient method,'' \emph{SIAM Journal on Control and Optimization}, vol.~59, no.~5, pp. 3887--3911, 2021.

\bibitem{bhandari_global_2022}
\BIBentryALTinterwordspacing
J.~Bhandari and D.~Russo, ``Global optimality guarantees for policy gradient methods,'' \emph{Operations Research}, vol.~72, no.~5, pp. 1906--1927, 2024. [Online]. Available: \url{https://doi.org/10.1287/opre.2021.0014}
\BIBentrySTDinterwordspacing

\bibitem{eftekhari_training_2020}
A.~Eftekhari, ``Training linear neural networks: Non-local convergence and complexity results,'' in \emph{International Conference on Machine Learning}.\hskip 1em plus 0.5em minus 0.4em\relax PMLR, 2020, pp. 2836--2847.

\bibitem{boyd2004convex}
S.~P. Boyd and L.~Vandenberghe, \emph{Convex Optimization}.\hskip 1em plus 0.5em minus 0.4em\relax Cambridge, England: Cambridge University Press, 2004.

\bibitem{chvatal1983linear}
V.~Chv{\'a}tal, \emph{Linear programming}.\hskip 1em plus 0.5em minus 0.4em\relax London, UK: Macmillan, 1983.

\bibitem{polyak1963gradient}
B.~T. Polyak, ``Gradient methods for minimizing functionals,'' \emph{USSR Computational Mathematics and Mathematical Physics}, vol.~3, no.~4, pp. 643--653, 1963.

\bibitem{karimi2016linear}
H.~Karimi, J.~Nutini, and M.~Schmidt, ``Linear convergence of gradient and proximal-gradient methods under the {P}olyak-{{\L}}ojasiewicz condition,'' in \emph{Machine Learning and Knowledge Discovery in Databases}, P.~Frasconi, N.~Landwehr, G.~Manco, and J.~Vreeken, Eds.\hskip 1em plus 0.5em minus 0.4em\relax Cham: Springer International Publishing, 2016, pp. 795--811.

\bibitem{watanabe2025revisiting}
Y.~Watanabe and Y.~Zheng, ``Revisiting strong duality, hidden convexity, and gradient dominance in the linear quadratic regulator,'' \emph{arXiv preprint arXiv:2503.10964}, 2025.

\bibitem{fazel2018global}
M.~Fazel, R.~Ge, S.~Kakade, and M.~Mesbahi, ``Global convergence of policy gradient methods for the linear quadratic regulator,'' in \emph{International Conference on Machine Learning}.\hskip 1em plus 0.5em minus 0.4em\relax PMLR, 2018, pp. 1467--1476.

\bibitem{sun2021learning}
Y.~Sun and M.~Fazel, ``Learning optimal controllers by policy gradient: {Global} optimality via convex parameterization,'' in \emph{2021 60th IEEE Conference on Decision and Control (CDC)}, 2021, pp. 4576--4581.

\bibitem{hu2023toward}
B.~Hu, K.~Zhang, N.~Li, M.~Mesbahi, M.~Fazel, and T.~Ba{\c{s}}ar, ``Toward a theoretical foundation of policy optimization for learning control policies,'' \emph{Annual Review of Control, Robotics, and Autonomous Systems}, vol.~6, no.~1, pp. 123--158, 2023.

\bibitem{mohammadi_convergence_2022}
\BIBentryALTinterwordspacing
H.~Mohammadi, A.~Zare, M.~Soltanolkotabi, and M.~R. Jovanovic, ``\BIBforeignlanguage{en}{Convergence and sample complexity of gradient methods for the model-free linear–quadratic regulator problem},'' \emph{\BIBforeignlanguage{en}{IEEE Transactions on Automatic Control}}, vol.~67, no.~5, pp. 2435--2450, May 2022. [Online]. Available: \url{https://ieeexplore.ieee.org/document/9448427/}
\BIBentrySTDinterwordspacing

\bibitem{mohammadi2021lack}
H.~Mohammadi, M.~Soltanolkotabi, and M.~R. Jovanovi{\'c}, ``On the lack of gradient domination for linear quadratic {G}aussian problems with incomplete state information,'' in \emph{2021 60th IEEE Conference on Decision and Control (CDC)}, 2021, pp. 1120--1124.

\bibitem{lojasiewicz1984gradients}
S.~{\L}ojasiewicz, ``Sur les trajectoires du gradient d'une fonction analytique. ({Trajectories} of the gradient of an analytic function),'' \emph{Semin. Geom., Univ. Studi Bologna}, vol. 1982/1983, pp. 115--117, 1984.

\bibitem{angeli2000characterization}
D.~Angeli, E.~D. Sontag, and Y.~Wang, ``A characterization of integral input-to-state stability,'' \emph{IEEE Transactions on Automatic Control}, vol.~45, no.~6, pp. 1082--1097, 2000.

\bibitem{sontag1989smooth}
E.~D. Sontag \emph{et~al.}, ``Smooth stabilization implies coprime factorization,'' \emph{IEEE Transactions on Automatic Control}, vol.~34, no.~4, pp. 435--443, 1989.

\bibitem{rautert_computational_1997}
\BIBentryALTinterwordspacing
T.~Rautert and E.~W. Sachs, ``\BIBforeignlanguage{en}{Computational design of optimal output feedback controllers},'' \emph{\BIBforeignlanguage{en}{SIAM Journal on Optimization}}, vol.~7, no.~3, pp. 837--852, Aug. 1997. [Online]. Available: \url{http://epubs.siam.org/doi/10.1137/S1052623495290441}
\BIBentrySTDinterwordspacing

\bibitem{hassan2021proximal}
S.~Hassan-Moghaddam and M.~R. Jovanovi{\'c}, ``Proximal gradient flow and {D}ouglas--{R}achford splitting dynamics: {G}lobal exponential stability via integral quadratic constraints,'' \emph{Automatica}, vol. 123, p. 109311, 2021.

\bibitem{gokhale2024proximal}
A.~Gokhale, A.~Davydov, and F.~Bullo, ``Proximal gradient dynamics: {M}onotonicity, exponential convergence, and applications,'' \emph{IEEE Control Systems Letters}, vol.~8, pp. 2853--2858, 2024.

\bibitem{sontag2013mathematical}
E.~D. Sontag, \emph{Mathematical control theory: deterministic finite dimensional systems}.\hskip 1em plus 0.5em minus 0.4em\relax Springer Science \& Business Media, 2013, vol.~6.

\end{thebibliography}
\end{spacing}

%%%%%%%%%%%%%%%%%%%%%%%%%%%%%%%%%%%%%%%%%%%%%%%%%%%%%%%%%%%%%%%%%%%%%%%%%%%%%%%%

\section*{APPENDIX}

\subsection{Proofs of Results from Section \ref{sec:convopt}}

\subsubsection*{Proof of Lemma \ref{lem:precompact}}
    Since $\cost$ is proper, then all level-sets of $\cost$ are compact, furthremore, since we adopt a gradient flow dynamics for the parameters $\x$ we have that
    \begin{equation*}
        \frac{d}{dt}\cost[\x] = \ip{\gradCost[\x],\dot x} = -\|\gradCost[\x]\|^2
    \end{equation*}
    which implies that the cost function is non-increasing along a trajectory of the gradient flow, which means that a solution $\sol[t,\x_0]$ is trapped inside the level-set $\cost[\x_0]$, which is compact and, therefore, implies that $\sol[t,\x_0]$ is pre-compact. \qed

\subsubsection*{Proof of Lemma \ref{lem:Conv2Min}} The proof of this lemma is a direct consequence of the results in Appendix A of \cite{de2024remarks,de2024convergence}, and thus it is omitted here for the brevity of the paper. %\textcolor{red}{Should we reproduce the full proof here ``for completeness'', since it is arXiv?}

\subsubsection*{Proof of Lemma \ref{lem:expconvPLI}}
    We first show that $\cost$ being $\mu$-\gPLI implies its gradient flow solution is $\mu$-GES. For that, let $\cost$ be a $\mu$-\gPLI cost function for some $\mu>0$, and define $g(\x) = \cost(x)-\minCost$, which implies that $\nabla g(\x) = \gradCost[\x]$. Then, compute
    \begin{equation*}
        \frac{d}{dt}g(\x) = -\|\nabla g(\x)\|^2\leq -\mu(g(\x))
    \end{equation*}
    which can be solved for the solution
    \begin{align*}
        g(x(t)) &\leq g(x(0))\mbox{e}^{-\mu t} \\
        \cost[\x(t)]-\minCost&\leq (\cost[\x_0]-\minCost)\mbox{e}^{-\mu t}
    \end{align*}
    proving the direct implication. To prove the converse, \ie that the gradient flow being $\mu$-GES implies that $\cost$ must be $\mu$-\gPLI, assume that is not the case. Then, notice that for all $\x_0\in\SSp$
    \begin{align*}
        -\|\nabla g(\x_0)\|^2&=\frac{\mathrm{d}}{\mathrm{d}t}g(\sol[t,\x_0])\bigg|_{t=0} \\ &= \lim_{t\rightarrow0}\frac{1}{t}\left[g(\sol[t,\x_0])-g(\x_0)\right] \\
        &\leq \lim_{t\rightarrow0}\frac{1}{t}\left[\mathrm{e}^{-\mu t}(g(\x_0)-g(\x_0)\right] \\
        &= g(\x_0)\lim_{t\rightarrow0}\frac{1}{t}\left[\mathrm{e}^{-\mu t}-1\right] \\
        &= -\mu g(\x_0),
    \end{align*}
    which implies that for all $\x_0\in\SSp$ it holds that
    \begin{align*}
        \sqrt{\mu g(x_0)} &\leq \|\nabla g(\x_0)\| \\
        \sqrt{\mu (\cost[x_0]-\minCost)} &\leq \|\gradCost[\x_0]\|
    \end{align*}
    reaching contradiction and finishing the proof. \qed

\subsubsection*{Proof of Lemma \ref{lem:PLIvssgPLI}}
    That Assumption \ref{def:PLI} implies Assumption \ref{def:sgPLI} follows immediately, just pick $\mu_\e=\mu$ for all $\e>0$. To show that if $\cost$ satisfies a \sgPLI but not a P\L I then the constant $\overline\mu_\e$ must go to zero for some sequence $\{\e_i\}$, $i=[1,2,\dots]$ with $\e_i>0$, assume that it is not the case, then there must exist some $\underline\mu>0$ such that for any sequence $\{\e_i\}$, $\lim_{i\rightarrow\infty}\overline\mu_{\e_i}\geq \underline\mu$, however this would mean that for all $\e>0$, $\mu_\e\geq \underline\mu$, which would mean that $\cost$ would satisfy Assumption \ref{def:PLI} with $\mu=\underline\mu$, reaching contradiction. \qed

\subsubsection*{Proof of Lemma \ref{lem:lpli-ksat}}

    The implication that 1)$\Rightarrow$2) is immediate since any class-$\mathcal{K}_{\mathrm{SAT}}$ is positive definite, and for a given $ar/(b+r)$ and $\e$, one can always find $\mu>0$ such that $ar/(b+r)\geq\mu r$ for $r\leq \e$. 
    
    To show that 2)$\Rightarrow$1), first notice that $\alpha$ being positive-definite and with a nonzero $\liminf$ implies that there exists a $\underline\ell>0$ such that $\|\gradCost[\x]\|\geq \underline\ell$ for all $x$ outside of $\SSp_\e$, the sublevelset in which the $\lPLI$ holds. For a class-$\mathcal{K}_{\mathrm{SAT}}$ function to lower-bound the gradient, it is enough for it to lower-bound the $\lPLI$ comparison function inside $\SSp_\e$ and to be smaller than $\underline \ell$ outside of it.

    For the first part, let $\mu$ be such that 
    \begin{equation*}
        \|\gradCost[\x]\|\geq \sqrt{\mu(\cost[\x]-\minCost)}
    \end{equation*}
    for all $\x\in\SSp_\e$. Then, any clas $\mathcal{K}_{\mathrm{SAT}}$ function that lower-bounds it in $\SSp_\e$ must satisfy
    \begin{align*}
        &\frac{ar}{b+r}<\mu r & \forall r\leq \e \\
        &a<\mu b + \mu r & \forall r\leq \e \\
        &r>\frac{a-\mu b}{\mu} & \forall r\leq \e
    \end{align*}
    which holds if and only if $a<\mu b$, since $r>0$ from the definition of the comparison function. Furthermore, notice that for the comparison function to be smaller than $\underline \ell$ outside of $\SSp_\e$, it is enough to have that
    \begin{align*}
        &\lim_{r\rightarrow\infty}\sqrt{\frac{ar}{b+r}} = \sqrt{a}\leq\ell & \forall r\geq \e \\
    \end{align*}
    this proves that there always exists $a,b>0$ such that $\cost$ satisfies \eqref{eq:classKPLI} with $\alpha(r) = ar/(b+r)$, completing the proof. \qed

\subsubsection*{Proof of Lemma \ref{lem:linexp-bnd}}
    This proof is done in two parts, first for the linear upper bound, and then for the exponential upper bound.
    
    Pick $\e>0$, $x_0\in\SSp$, and $\underline t$ as described in the lemma's statement. Furthermore, let $\underline c$ be the smallest $c>0$ such that $\gradCost[\x]\leq c$ for all $x\in\SSp$. Then, notice that for any $t>0$:
    \begin{align*}
        \dot\cost(\phi( t,x_0))=&-\|\gradCost[ t,x_0)]\|^2 \\ \ge&-\underline c^2
    \end{align*}
    which implies that
    \begin{align*}
        \int_{t_1}^{t_2}\dot\cost(\phi(s,x_0))\mathrm{d} s \geq& -\underline c^2 (t_2-t_1) \\
        \cost[\sol(t_2,x_0)]-\cost[\sol(t_1,x_0)]\geq& -\underline c^2 (t_2-t_1).
    \end{align*}
    Then, pick $t_2=\underline t$ and $t_1=t$ for any $t<\underline t$, we have that
    \begin{align*}
        &\cost[\sol(\underline t,x_0)]-\cost[\sol(t,x_0)]\geq -\underline c^2 (\underline t-t) \\
        &\cost[\sol(t,x_0)]-\minCost\leq (\cost[\sol(\underline t,x_0)]-\minCost)-\underline c^2 (t-\underline t)
    \end{align*}
    which recovers the first statement of the lemma with $\delta=\underline c^2$ and $\delta_\e=\underline t\underline c^2+\e+\minCost$.

    For the second statement, consider that $\cost$ satisfies \eqref{eq:classKPLI} with a class $\mathcal{K}_\mathrm{SAT}$ comparison function as in \eqref{eq:Ksatdef}. Let $r(x) = \cost[\x]-\minCost$ and notice that for all $\x\in\SSp_\e$
    %\
    \begin{align*}
        \|\gradCost[x]\|^2\geq& \frac{a}{b+r(x)}r(x) \\ \geq & \frac{a}{b+\e}(\cost[\x]-\minCost) \\ =& \mu_\e(\cost[\x]-\minCost)
    \end{align*}
    Then, applying Lemma \ref{lem:expconvPLI} for $\SSp=\SSp_\e$ and $\mu=\mu_\e$ recovers the second statement of the lemma, completing the proof. \qed

\subsection{Proofs of Results from Section \ref{sec:LQRPO}}

\subsubsection*{Proof of Lemma \ref{lem:hgd-exists}}

    Let $(A,B)$ be a controllable pair, then we want to design a family of stabilizing feedback gain matrices $u=-\Kr x$,
    parametrized by $\rho\in(0,\infty)$, such that
    $\|\Kr\| \rightarrow+\infty$ as $\rho\rightarrow+\infty$.
    
    Consider first the single-input case.
    We build $\Kr$ as the unique feedback that assigns the following
    characteristic polynomial: 
    \begin{equation*}
        p(s) = (s + \rho)^n
    \end{equation*}
    Expanding this polynomial using the binomial theorem:
    \begin{equation*}
        p(s) = s^n + \binom{n}{1} \rho s^{n-1} + \binom{n}{2} \rho^2 s^{n-2} + \dots + \binom{n}{n} \rho^n.
    \end{equation*}
    Thus, the coefficients of the characteristic polynomial are:
    \begin{equation*}
        \alpha_i = \binom{n}{i} \rho^i, \quad \text{for } i = 1,2, \dots, n.
    \end{equation*}
    
    Ackermann's formula for the state feedback gain is given
    by\footnote{to check that I am not missing a negative sign somewhere!}:
    \begin{equation*}
        \Kr = \begin{bmatrix} 0 & \cdots & 0 & 1 \end{bmatrix} P^{-1} \Phi(A),
    \end{equation*}
    where
    \[
    P = \begin{bmatrix} B & AB & A^2B & \dots & A^{n-1}B \end{bmatrix}
    \]
    is the controllability matrix,
    Now, $\Phi(A)$ is obtained by substituting $A$ into the characteristic polynomial:
    \begin{equation*}
        \Phi(A) = A^n + \binom{n}{1} \rho A^{n-1} + \binom{n}{2} \rho^2 A^{n-2} + \dots + \binom{n}{n} \rho^n I.
    \end{equation*}
    
    Thus, the state feedback gain matrix $\Kr$ is:
    \begin{align*}
        \Kr =& \begin{bmatrix} 0 & \cdots & 0 & 1 \end{bmatrix} P^{-1} \\&\times\left( A^n + \binom{n}{1} \rho A^{n-1} + \dots + \binom{n}{n} \rho^n I \right).
    \end{align*}
    This ensures that the closed-loop system has the desired eigenvalues at $-\rho$ with multiplicity $n$, leading to an exponentially stable system with the decay rate $e^{-\rho t}$, and guaranteeing that the eigenvalues of the closed-loop system are bounded away from the imaginary axis.
    
    Notice that this expression shows that $\|\Kr\| \rightarrow+\infty$
    as $\rho\rightarrow+\infty$, which from coerciveness of the cost $J(\cdot)$ implies that $\lim_{\rho\rightarrow\infty}J(K^*+\Kr)=\infty$.
    
    In the case $m>1$ and $(A,B)$ controllable, we know from the proof of
    Theorem 13 in \cite{sontag2013mathematical} that there is a matrix $F$ and a vector $v$ such that $(A-BF,Bv)$ is controllable. We apply Ackermann's formula to this single-input system, obtaining a matrix family $K(\rho)$ as above. Now $\Kr := F + vK(\rho)$ is as desired. \qed

\subsubsection*{Proof of Lemma \ref{lem:BndGradJ}}

    Let $\Kr$ be any high gain curve of $\Kstbl$ as defined in Definition \ref{def:hgd}, such that $B\Ke$ is diagonalizable. From coerciveness of the cost $J(\cdot)$ and the fact that $\Kr$ is bounded away from the imaginary axis, we can conclude that $\lim_{\rho\rightarrow\infty}\|\Kr\| =\infty$. From this, let $\rt:\re_+\rightarrow\re_{++}$ be any function such that 
    $\lim_{\rho\rightarrow\infty}\frac{\|\Kr\|}{\rt(\rho)} <\infty$ (we will write $\rt$ instead of $\rt(\rho)$ for simplicity). 
    
    We will assume in this proof that all unbounded eigenvalues of $B\Kr$ grow to infinity at the same rate as $\rt$. The proof can be extended for the case where the eigenvalues of $B\Kr$ grow to infinity at different rates, however, it would grow significantly more complicated. 
    
    Furthermore, notice the high gain curve constructed in the proof of \ref{lem:hgd-exists} satisfies this extra assumption, and so it is enough for the goals of the paper. More than that, the high gain curve constructed there is such that only one eigenvalue of $B\Kr$ is going to infinity, with all others remaining constant (equivalent to $p=1$ in the next equation). 
    
    Then, for any $\rho>0$ write
    \begin{equation*}
        B\Kr = V\begin{bmatrix}
            \rt\Lambda_{\rho} & ~~0 \\ 0 & ~~\rt Z_{\rho}
        \end{bmatrix}V_{\rho}^{-1} = \rt V_{\rho}\bar\Lambda_{\rho} V_{\rho}^{-1},
    \end{equation*}
    where $Z_{\rho}\in\re[q\times q]$ is a diagonal matrix such that $\lim_{\rho\rightarrow\infty}\|\rt Z_{\rho}\|<\infty$, and matrix $\Lambda_{\rho}\in\mathbb{R}^{\p\times \p}$, is 
    a diagonal matrix such that $\lim_{\rho\rightarrow\infty}\sigma_{\min}(\Lambda_\rho)>0$. In other words, $\rt Z_\rho$ collects the eigenvalues of $B\Kr$ that either are zero or constant as $\rho\rightarrow\infty$, while $\rt\Lambda_\rho$ collects only the eigenvalues that grow to infinity at the same rate as $\rt$ as $\rho\rightarrow\infty$.
    
    Next, let $A^*=A-BK^*$, and consider the Lyapunov equation for $Y_\rho := Y_{\Kr}$ and write
    \begin{align*}
        (A^*-B\Kr)Y_{\rho}+Y_{\rho}(A^*-B\Kr)^\top+I_n&=0 \\
        (\Abs-\rt\bar\Lambda_{\rho})\Ykb+\Ykb(\Abs-\rt\bar\Lambda_{\rho})^\top +\I &= 0
    \end{align*}
    where $\Abs = V_{\rho}^{-1}A^*V_{\rho}$, $\Ykb = V_{\rho}^{-1}Y_{\rho}V_{\rho}^{-\top}$, and $\I = V_{\rho}^{-1}V_{\rho}^{-\top}$. From this, write
    
    \begin{align*}
        &\left(\begin{bmatrix}
            \Abs_{11} & \Abs_{12} \\ \Abs_{21} & \Abs_{22}-\rt Z_{\rho}
        \end{bmatrix}-\begin{bmatrix}
            \rt\Lambda_{\rho} & ~~0 \\ 0 & ~~0
        \end{bmatrix}\right)\begin{bmatrix}
            \Ykb[1] & \Ykb[2] \\ \Ykb[2]^\top & \Ykb[3]
        \end{bmatrix} \\&+ \begin{bmatrix}
            \Ykb[1] & \Ykb[2] \\ \Ykb[2]^\top & \Ykb[3]
        \end{bmatrix}\left(\begin{bmatrix}
            \Abs_{11} & \Abs_{12} \\ \Abs_{21} & \Abs_{22}-\rt Z_{\rho}
        \end{bmatrix}-\begin{bmatrix}
            \rt\Lambda_{\rho} & ~~0 \\ 0 & ~~0
        \end{bmatrix}\right)^\top \\&+ \begin{bmatrix}
            \I[1] & \I[2] \\ \I[2]^\top & \I[3]
        \end{bmatrix} = 0
    \end{align*}
    with $\Abs_{11}\in\re[p\times p]$, $\Abs_{12}\in\re[p\times q]$, $\Abs_{21}\in\re[q\times p]$ and $\Abs_{22}\in\re[q\times q]$, and similar partitions for $\Ykb$ and $\I$. Notice that the eigenvalues of $\Abs_{22}-\rt Z_{\rho}$ are bounded away from zero (since $\Kr\in\Kstbl$) and bounded above, implying that $\lim_{\rho\rightarrow\infty}\|\Abs_{22}-\rt Z_{\rho}\|<\infty$. Therefore, for convenience, we will override the notation as $\Abs_{22} = \Abs_{22}-\rt Z_{\rho}$ for the remainder of the proof. 
    
    Then, opening this results in the following linear matrix equalities
    \begin{align}
        (\Abs_{11}-\rt\Lambda_{\rho})\Ykb[1]+\Abs_{12}\Ykb[2]^\top\nonumber+\Ykb[1](\Abs_{11}-\rt\Lambda_{\rho})^\top&\\+\Ykb[2] (\Abs_{12})^\top + \I[1]&=0\label{eq:lmi1}\\
        \Ykb[2]^\top(\Abs_{11}-\rt\Lambda_{\rho})^\top+\Ykb[3](\Abs_{12})^\top\nonumber& \\+\Abs_{21}\Ykb[1]+\Abs_{22}\Ykb[2]^\top + \I[2]^\top &= 0 \label{eq:lmi21}\\
        (\Abs_{11}-\rt\Lambda_{\rho})\Ykb[2]+\Abs_{12}\Ykb[3]+\Ykb[1](\Abs_{21})^\top\nonumber&\\+\Ykb[2]
        (\Abs_{22})^\top + \I[2]&=0\label{eq:lmi12}\\
        \Abs_{21}\Ykb[2]+\Ykb[2]^\top(\Abs_{21})^\top +\Abs_{22}\Ykb[3] \nonumber&\\+ \Ykb[3](\Abs_{22})^\top + \I[3] &= 0 \label{eq:lmi3}
    \end{align}
    
    For the next steps of the derivation, we will vectorize the equations and use the following two properties of the vectorization operator, $$\vect{A+B} = \vect{A} + \vect{B},$$ and $$\vect{ABC} = (C^\top \otimes A)\vect{B}.$$ Furthermore, we denote by $\Pm[k,l]$ the permutation matrix of appropriate dimensions that satisfies the following identity $$\vect{A^\top} = \Pm[k,l]\vect{A}$$ for $A\in\re[k\times l]$. We will also denote $\vA[k]_{ij}:=I_k\otimes \Abs_{ij}+(\Abs_{ij}\otimes I_k)\Pm$ (where the dimension of $\Pm$ is implicit by context) and $\vL[k] := I_k\otimes \Lambda + \Lambda\otimes I_k$. With these stabilished, vectorize \eqref{eq:lmi1} and \eqref{eq:lmi3} to obtain
    \begin{align*}
        \vect{\Ykb[1]} &= -(\vA[p]_{11}-\rt\vL[p])^{-1}(\vA[p]_{12}\Pm[p,q]\vect{\Ykb[2]}\\&~~~~~~~+\vect{\I[1]})\\
        \vect{\Ykb[3]} &= -(\vA[q]_{22})^{-1}(\vA[q]_{21}\vect{\Ykb[2]}+\vect{\I[3]}).
    \end{align*}
    We next vectorize \eqref{eq:lmi12} to obtain
    \begin{align*}
        &\left(I_q\otimes\Abs_{11}+\Abs_{22}\otimes I_p-\rt(I_q\otimes\Lambda)\right)\vect{\Ykb[2]} \\&= -(I_q\otimes\Abs_{12})\vect{\Ykb[3]}-(\Abs_{21}\otimes I_p)\vect{\Ykb[1]}-\vect{\I[2]} \\
        &(\overline\vA_{11}+\overline\vA_{22}-\rt\overline\vL)\vect{\Ykb[2]} \\&= -\overline\vA_{12}\vect{\Ykb[3]}-\overline\vA_{21}\vect{\Ykb[1]}-\vect{\I[2]}
    \end{align*}
    %
    % \begin{align*}
    %     (\overline\vA_{11}+\overline\vA_{22}-\rho\overline\vL)\vect{\Ykb[2]} = -\vA[q]_{12}\vect{\Ykb[3]}-\vA[p]_{21}\vect{\Ykb[1]}-(I+\Pm)\vect{\I[2]}
    % \end{align*}
    %
    % where $$\overline\vA_{11} = (\Abs_{11}\otimes I_q) \Pm[p,q]+(I_q\otimes \Abs_{11}),$$ $$\overline \vA_{22} = (I_p\otimes \Abs_{22})\Pm[p,q]+(\Abs_{22}\otimes I_p),$$ and $$\overline\vL = (\Lambda\otimes I_q)\Pm[p,q]+(I_q\otimes\Lambda).$$
    
    We then use the two identities for $\vect{\Ykb[1]}$ and $\vect{\Ykb[3]}$ and solve for $\vect{\Ykb[2]}$ as
    \begin{align*}
        \vect{\Ykb[2]} =& \big[\overline \vA_{11}+\overline\vA_{22}-\rt\overline\vL-\overline\vA_{12}(\vA[q]_{22})^{-1}\vA[q]_{21}\\&-\overline\vA_{21}(\vA[p]_{11}-\rt\vL[p])^{-1}\vA[p]_{12}\Pm[p,q]\big]^{-1} \\ &\times \big(\overline\vA_{12}(\vA[q]_{22})^{-1}\vect{\I[3]}\\&+\overline\vA_{21}(\vA[p]_{11}-\rt\vL[p])^{-1}\vect{\I[1]}-\vect{\I[2]}\big)
    \end{align*}
    Notice from the expression above that $$\lim_{\rho\rightarrow\infty} \vect{\Ykb[2]} = 0,$$ which when applied to \eqref{eq:lmi1} implies that $$\lim_{\rho\rightarrow\infty}\vect{\Ykb[1]}=0.$$ However, notice also that 
    \begin{align*}
        \lim_{\rho\rightarrow\infty} \rt\vect{\Ykb[2]} &= (\overline\vL)^{-1}\big(\overline\vA_{12}(\vA[q]_{22})^{-1}\vect{\I[3]}\\&~~~~~~~~~~~~~-\vect{\I[2]}\big)\\&=:\vect{\Ykb[2]^\infty},
    \end{align*}
    and that $$\lim_{\rho\rightarrow\infty}\rt\Ykb[1]=:\Ykb[1]^\infty = (\vL[p])^{-1}(\vect{\I[1]}).$$ This means that both terms $\Ykb[1]$ and $\Ykb[2]$ are going to zero at a rate of $1/\rt$.
    %\LC{Can we ensure that $\overline\vL$ and $\vL{^p}$ are invertible?}
    
    Next, define $\Ykb[3]^\infty:=\lim_{\rho\rightarrow\infty}\Ykb[3]$ (without the $\rt$ this time), and notice from \eqref{eq:lmi3} and from the fact that $\Ykb[2]$ goes to $0$ as $\rho\rightarrow\infty$, that it must solve the following LME
    \begin{equation*}
        \Abs_{22}\Ykb[3]^\infty+\Ykb[3]^\infty(\Abs_{22})^\top+\I[3]=0,
    \end{equation*}
    Therefore, we finally conclude that
    \begin{equation*}
        \lim_{\rho\rightarrow\infty}\Ykb = \begin{bmatrix}
            0 & 0 \\ 0 & ~~\Ykb[3]^\infty
        \end{bmatrix}.
    \end{equation*}
    
    Next, consider the Lyapunov equation for $\PrKe:=P_{\Kr}$ %, and \textcolor{red}{by a very similar argument}\footnote{I have it in my notes, I just have not transcribed it yet, will do as soon as I wake up.}
    \begin{equation*}
        (A^*-B\Kr)^\top \PrKe+\PrKe(A^*-B\Kr)+\Kr^\top R\Kr=0,
    \end{equation*}
    define $\Pkb=V^\top \PrKe V$ and notice that $(1/\rt^2)V^\top \Ke^\top R \Ke V$ has the following block structure
    \begin{equation*}
        \frac{1}{\rt^2}V_{\rho}^\top \Kr^\top R \Kr V_{\rho}=\begin{bmatrix}
            \Rp & ~~0 \\ 0 & ~~\Op
        \end{bmatrix}
    \end{equation*}
    where for all $\rho$, $\Rp\in\PD[p]$ and is bounded, and $\lim_{\rho\rightarrow\infty}\rt^2\Op=\Op^\infty$ with $\|\Op^\infty\|<\infty$. To see this, let $V_{\rho} = [V_{\rho,\p},~V_{\rho,\q}]$ be the first $\p$ and last $\q$ columns of $V_{\rho}$. Notice that necessarily, $\lim_{\rho\rightarrow\infty}(1/\rt)\Kr V_{\rho,\q}=0$, since $B\Kr V_{\rho,\q}=V_{\rho,q}\rt Z_{\rho}$ and $B$ is full rank. The block structure follows from that. From here, we rewrite the Lyapunov equation as
    \begin{align*}
        &\left(\begin{bmatrix}
            \Abs_{11} & \Abs_{12} \\ \Abs_{21} & \Abs_{22}
        \end{bmatrix}-\begin{bmatrix}
            \rt\Lambda & ~~0 \\ 0 & ~~0
        \end{bmatrix}\right)^\top\begin{bmatrix}
            \Pkb[1] & \Pkb[2] \\ \Pkb[2]^\top & \Pkb[3]
        \end{bmatrix} \\&+ \begin{bmatrix}
            \Pkb[1] & \Pkb[2] \\ \Pkb[2]^\top & \Pkb[3]
        \end{bmatrix}\left(\begin{bmatrix}
            \Abs_{11} & \Abs_{12} \\ \Abs_{21} & \Abs_{22}
        \end{bmatrix}-\begin{bmatrix}
            \rt\Lambda & ~~0 \\ 0 & ~~0
        \end{bmatrix}\right) \\&+ \begin{bmatrix}
            \rt^2\Rp & ~~0 \\ 0 & ~~\rt^2\Op
        \end{bmatrix} = 0
    \end{align*}
    
    From this, we open the equation above into the following LMEs
    \begin{align}
        (\Abs_{11}-\rt\Lambda)^\top\Pkb[1]+(\Abs_{21})^\top\Pkb[2]^\top\nonumber&\\+\Pkb[1](\Abs_{11}-\rt\Lambda)+\Pkb[2] \Abs_{21} + \rt^2\Rp&=0\label{eq:lme1}\\
        \Pkb[2]^\top(\Abs_{11}-\rt\Lambda)+\Pkb[3]\Abs_{21}\nonumber&\\+(\Abs_{12})^\top\Pkb[1]+(\Abs_{22})^\top\Pkb[2]^\top &= 0 \label{eq:lme21}\\
        (\Abs_{11}-\rt\Lambda)^\top\Pkb[2]+(\Abs_{21})^\top\Pkb[3]\nonumber&\\+\Pkb[1]\Abs_{12}+\Pkb[2] \Abs_{22} &=0\label{eq:lme12}\\
        (\Abs_{12})^\top\Pkb[2]+\Pkb[2]^\top \Abs_{12} +(\Abs_{22})^\top\Pkb[3] \nonumber&\\+ \Pkb[3]\Abs_{22}+\rt^2\Op&= 0 \label{eq:lme3}
    \end{align}
    
    Then, similarly to the procedure for $\Ykb$, we vectorize \eqref{eq:lme1} and \eqref{eq:lme3} and solve for $\Pkb[1]$ and $\Pkb[3]$ as follows
    
    \begin{align*}
        \vect{\Pkb[1]} &= -(\vA[p\top]_{11}-\rt\vL[p])^{-1}(\vA[p\top]_{21}\Pm[p,q]\vect{\Pkb[2]}\\&~~~~~~~+\rt^2\vect{\Rp})\\
        \vect{\Pkb[3]} &= -(\vA[q\top]_{22})^{-1}(\vA[q\top]_{12}\vect{\Pkb[2]}+\rt^2\vect{\Op}).
    \end{align*}
    We then vectorize \eqref{eq:lme12} to obtain
    \begin{align*}
        &\left(I_q\otimes\Abs_{11}+\Abs_{22}\otimes I_p-\rt(I_q\otimes\Lambda)\right)^\top\vect{\Pkb[2]} \\&~~~= -(I_q\otimes\Abs_{21})^\top\vect{\Pkb[3]}-(\Abs_{12}\otimes I_p)^\top\vect{\Pkb[1]} \\
        &(\overline \vA^\top_{11}+\overline\vA^\top_{22}-\rt\overline\vL {^\top})\vect{\Pkb[2]} \\&~~~= -\overline\vA_{21}^\circ\vect{\Pkb[3]}-\overline\vA_{12}^\circ\vect{\Pkb[1]}
    \end{align*}
    %\LC{$\overline\vA_{21}^\top$ and $\overline\vA_{12}^\top$ are not the transpose of $\overline\vA_{21}$ and $\overline\vA_{12}$. We can use different notations.}
    %
    %where $\overline\vA_{12}^\circ:=(\Abs_{12}\otimes I_p)^\top$ and $\overline\vA_{21}^\circ:=(I_q\otimes\Abs_{21})^\top$. 
    We then use the two identities for $\vect{\Pkb[1]}$ and $\vect{\Pkb[3]}$ and solve for $\vect{\Pkb[2]}$ as
    \begin{align*}
        \vect{\Pkb[2]} =& \Big[\overline \vA_{11}^\top+\overline\vA_{22}^\top-\rt\overline\vL-\overline\vA_{21}^\circ(\vA[q\top]_{22})^{-1}\vA[q\top]_{12}\\&-\overline\vA_{12}^\circ(\vA[p\top]_{11}-\rt\vL[p\top])^{-1}\vA[p\top]_{21}\Pm[pq]\Big]^{-1} \\ &\times \Big(\overline\vA_{12}^\circ(\vA[p\top]_{11}-\rt\vL[p])^{-1}\vect{\Rp}\rt^2\\&+\overline\vA_{21}^\circ(\vA[q\top]_{22})^{-1}\vect{\Op}\rt^2\Big).
    \end{align*}
    
    From the expression above, notice that $$\lim_{\rho\rightarrow\infty}\vect{\Pkb[2]} = (\overline\vL)^{-1}\overline\vA_{21}^\circ(\vL[p])^{-1}\vect{\Rp}=:\vect{\Pkb[2]^\infty},$$ which in turn implies that $$\lim_{\rho\rightarrow\infty}\frac{1}{\rho}\vect{\Pkb[1]} = (\vL[p])^{-1}\vect{\Rp}=:\vect{\Pkb[1]^\infty},$$ and that 
    \begin{align*}
        \lim_{\rho\rightarrow\infty}\vect{\Pkb[3]} &= -(\vA[q\top]_{22})^{-1}(\vA[q\top]_{12}\vect{\Pkb[2]^\infty}+\vect{\Op^\infty})\\&=:\vect{\Pkb[3]^\infty}
    \end{align*}
    % Following a very similar derivation as done for $\Ykb$\footnote{I have it in my hand notes, I will transcribe them tomorrow after I sleep a little} we get that $\lim_{\rho\rightarrow\infty}(1/\rho)\Pkb[1]$ is finite and non-zero and $\lim_{\rho\rightarrow\infty}(1/\rho)\Pkb[2]=0$ and $\lim_{\rho\rightarrow\infty}(1/\rho)\Pkb[3]=0$. 
    With these results, consider the expression of the gradient
    {\small
    \begin{align*}
        \nabla J(\Kr) &= 2(R\Kr-B^\top \PrKe)Y_{\rho} \\
        \nabla J(\Kr)V^{-\top}&=2(R\Kr V - B^\top V^{-\top}V^\top\PrKe V)V^{-1}Y_{\rho}V^{-\top}\\
        \nabla J(\Kr)V^{-\top}&=2(R\Keb - \Bb(1/\rho)\Pkb)\rho\Ykb.
    \end{align*}}
    In an attempt to improve clarity, we will look at each term of $\nabla J(\Kr)$ separately. First define
    \begin{align*}
        \Kr V &= [\Kr V_p, ~~ \Kr V_q]\\
        &=:[\Kt_p(\rho), ~~\Kt_q(\rho)]
    \end{align*}
    and notice that $\lim_{\rho\rightarrow\infty}\|(1/\rt)\Kt_p\|<0$, and $\lim_{\rho\rightarrow\infty}\|\Kt_q\|<\infty,$
    %
    % \begin{align*}
    %     \lim_{\rho\rightarrow\infty}(1/\rt)\Kr V &= \lim_{\rho\rightarrow\infty} (1/\rt)\Kr[V_p ~~V_q] \\&= \lim_{\rho\rightarrow\infty}(1/\rt)[\Kr V_p ~~\Kr V_q]\\&=:[\Kt ~~0],
    % \end{align*}
    %
    which implies that
    {\footnotesize
    \begin{align*}
        &(1/\rt)\Keb(\rt\Ykb) \\&= \begin{bmatrix}
            ((1/\rt)\Kt_p(\rho))(\rt\Ykb[1])+\Kt_q\Ykb[2]^\top & ((1/\rt)\Kt_p)(\rt\Ykb[2])+\Kt_q\Ykb[3]
        \end{bmatrix}.
    \end{align*}}
    Then, let $\lim_{\rho\rightarrow\infty}(1/\rt)\Kt_p=:\Kt_p^\infty$, and $\lim_{\rho\rightarrow\infty}\Kt_q=:\Kt_q^\infty$, and notice that 
    \begin{align*}
        &\lim_{\rho\rightarrow\infty} (1/\rt)R\Keb(\rt\Ykb) \\&= R\begin{bmatrix}
            \Kt_p^\infty\Ykb[1]^\infty & ~~\Kt_p^\infty\Ykb[2]^\infty+\Kt_q^\infty\Ykb[3]^\infty
        \end{bmatrix}\in\re[m\times n].
    \end{align*}
    
    For the second term write
    \begin{align*}
        \Pkb \Ykb = \begin{bmatrix}
            \Pkb[1]\Ykb[1]+\Pkb[2]\Ykb[2]^\top & ~~\Pkb[1]\Ykb[2]+\Pkb[2]\Ykb[3] \\ \Pkb[2]^\top\Ykb[1]+\Pkb[3]\Ykb[2]^\top & ~~\Pkb[2]^\top\Ykb[2]+\Pkb[3]\Ykb[3]
        \end{bmatrix}
    \end{align*}
    which allows us to analyse each term individually. From the previously obtained result we can say that 
    \begin{align*}
        &\lim_{\rho\rightarrow\infty}\Pkb[1]\Ykb[1] = \Pkb[1]^\infty\Ykb[1]^\infty\\
        &\lim_{\rho\rightarrow\infty}\Pkb[2]\Ykb[2]^\top = 0\\
        &\lim_{\rho\rightarrow\infty}\Pkb[2]^\top\Ykb[2] = 0\\
        &\lim_{\rho\rightarrow\infty}\Pkb[1]\Ykb[2] = \Pkb[1]^\infty\Ykb[2]^\infty\\
        &\lim_{\rho\rightarrow\infty}\Pkb[2]\Ykb[3] = \Pkb[2]^\infty\Ykb[3]^\infty\\
        &\lim_{\rho\rightarrow\infty}\Pkb[2]^\top\Ykb[1] = 0\\
        &\lim_{\rho\rightarrow\infty}\Pkb[3]\Ykb[2]^\top = 0\\
        &\lim_{\rho\rightarrow\infty}\Pkb[3]\Ykb[3] = \Pkb[3]^\infty\Ykb[3]^\infty.
    \end{align*}
    
    As such, the result of the limit for this term is given by
    \begin{equation*}
        \lim_{\rho\rightarrow\infty}\bar B^\top\Pkb\Ykb = \begin{bmatrix}
            \Pkb[1]^\infty\Ykb[1]^\infty & ~~\Pkb[1]^\infty\Ykb[2]^\infty+\Pkb[2]^\infty\Ykb[3]^\infty \\ 0 & ~~\Pkb[3]^\infty\Ykb[3]^\infty
        \end{bmatrix}
    \end{equation*}
    
    Therefore, the limit of both terms of $$\lim_{\rho\rightarrow\infty}\nabla J(\rho\Ke)V^{-\top}$$ exist and are finite, which implies that the gradient itself exists and is finite independent of the chosen unbounded direction, completing the proof with the assumption that all unbounded eigenvalues of $B\Kr$ grow to infinity at the same rate as $\rt$.
    
    If there were two distinct sets of eigenvalues of $B\Kr$ going to infinity at two distinct rates $\rt_1$ and $\rt_2$, one can rederive this proof by writing the SVD of $B\Kr$ as
    \begin{equation*}
        B\Kr = V\begin{bmatrix}
            \rt_1\Lambda_{1,\rho} & ~~0 & ~~0 \\ 0 &\rt_2\Lambda_{2,\rho} & ~~ 0 \\ 0 & ~~0 & Z_{\rho}
        \end{bmatrix}V_{\rho}^{-1} = \rt V_{\rho}\bar\Lambda_{\rho} V_{\rho}^{-1},
    \end{equation*}
    and breaking the LMEs into $3\times 3$ blocks instead of $2\times 2$, and verifying that the product of the terms of $P_{\Kr}$ and $Y_{\Kr}$ will converge to either zero or a constant still. As mentioned before, this would result in many more terms in the analysis, complicating it significantly, so we refrain from providing this proof. \qed

\subsubsection*{Proof of Corollary \ref{cor:LQRnoPLI}}

    Let $\Kr$ be an unbounded direction and assume $J$ is \gPLI with constant $\mu>0$. Then for all $\rho>0$ one must have
    \begin{equation*}
        \|\nabla J[\Kr]\|\geq \sqrt{\mu(J(\Kr)-J^*)}.
    \end{equation*}
    However, $\lim_{\rho\rightarrow\infty}J(\Kr)-J^*=\infty$ by assumption, and $\lim\|\nabla J(\Kr)\|<\infty$ from Lemma \ref{lem:BndGradJ}, which reaches contradiction. \qed

\subsubsection*{Proof of Lemma \ref{lem:GradJUnbd}}

For a given $K\in\Kstbl$, let $(A-BK) = V \Lambda V^\top$, be an eigendecomposition of $(A-BK)$ with $\lm$ being its eigenvalue with the largest real part. Then, let $\{K_i\}$, $i=1,2,\dots$ be an infinite sequence of feedback matrices $K_i\in\Kstbl$ such that $\Re(\lm_i)>\Re(\lm_j)$ for any $i>j$ with all other eigenvalues remaining the same, and let $\lim_{i\rightarrow\infty}\Re(\lm_i)=0$. Then, let $v_i$ be such that $(A-BK_i)v_i = \Re(\lm_i)v_i$, and write the following
\begin{align*}
    &v_i^\top (P_K(A-BK_i)+(A-BK_i)^\top P_K\\&~~~~~+K_i^\top RK_i+Q)v_i=0\\
    &v_i^\top (P_K(A-BK_i)+(A-BK_i)^\top P_K)v_i\\&~~~~~~=-v_i^\top(K_i^\top RK_i+Q)v_i \\
    &2\Re(\lm_i)v_i^\top P_K v_i \leq -v_i^\top Q v_i \\
    &\trace{P_i}\geq -\frac{\underline\lambda_Q}{2\Re(\lm_i)}
\end{align*}
Which, at the limit, goes to infinity. \qed

\subsubsection*{Proof of Lemma \ref{lem:LQRClassKPLIub}}
Notice that any finite point $K$ in $\partial\Kstbl_\delta := \{K\in\re[m\times n]~|~\Re(\lambda_{\max}(A-BK))=-\delta\}$ is strictly inside $\Kstbl$ and therefore has a finite value for $\nabla J(K)$. Furthermore, any trajectory $\Kr$ such that $\lim_{\rho\rightarrow\infty}\Kr\in\Kstbl_\delta$ and that $\lim_{\rho\rightarrow\infty}J(\Kr)-J^* = \infty$ is a high gain curve as per Definition \ref{def:hgd} and thus converges to a finite gradient.

%%%%%%%%%%%%%%%%%%%%%%%%%%%%%%%%%%%%%%%%%%%%%%%%%%%%%%%%%%%%%%%%%%%%%%%%%%%%%%%%

\end{document}